\documentclass[11pt]{article}

\usepackage{amssymb,amsmath,amsfonts,amsthm}
\usepackage[dvipsnames]{xcolor}

\usepackage{caption2}

\renewcommand{\subsection}{\subsubsection}

\renewcommand\appendix{\par
  \setcounter{section}{0}
   \renewcommand\thesection{Appendix \Alph{section}.}
 }

\newtheorem{theorem}{Theorem}[section]

\newtheorem{lemma}{Lemma}[section]
\newtheorem{proposition}{Proposition}[section]

\newtheorem{remark}{Remark}[section]

\def\calE{\mathcal E }
\def\calH{\mathcal H }

\newcommand{\nt}{|\hspace{-0.3pt}|\hspace{-0.3pt}|}

\begin{document}
\title{\bf Stability of an incompressible plasma-vacuum interface  with displacement current in vacuum}

\author{
Alessandro Morando$^{\ast}$\;,
Paolo Secchi$^{\dagger}$\;,
Yuri Trakhinin$^{\ddagger}$\;,
Paola Trebeschi$^{\S}$\\\\
\small{$\ast,\dagger,\S$ DICATAM, Sezione di Matematica, Universit\`a di Brescia} \\ \small{Via Valotti 9, 25133 Brescia, Italy}\\
\small{e-mail: alessandro.morando@unibs.it, paolo.secchi@unibs.it, paola.trebeschi@unibs.it}\\\\
\small{$\ddagger$ Sobolev Institute of Mathematics}\\ \small{Koptyug av. 4, 630090 Novosibirsk, Russia}\\
\centerline{\small{and}}\\
\small{Novosibirsk State University}\\\small{Pirogova str. 1, 630090 Novosibirsk, Russia}\\
\small{e-mail: trakhin@math.nsc.ru}}

\date{ }

%\author{{\bf Alessandro Morando}\\
%DICATAM, Sezione di Matematica, Universit\`a di Brescia \\ Via Valotti, 9, 25133 Brescia, Italy\\
%E-mail: alessandro.morando@unibs.it
%\and
%{\bf Paolo Secchi}\\
%DICATAM, Sezione di Matematica, Universit\`a di Brescia \\ Via Valotti, 9, 25133 Brescia, Italy\\
%E-mail: paolo.secchi@unibs.it
%\and
%{\bf Yuri Trakhinin}\\
%Sobolev Institute of Mathematics, Koptyug av. 4, 630090 Novosibirsk, Russia\\
%and\\
%Novosibirsk State University, Pirogova str. 1, 630090 Novosibirsk, Russia\\
%E-mail: trakhin@math.nsc.ru
%\and
%{\bf Paola Trebeschi}\\
%DICATAM, Sezione di Matematica, Universit\`a di Brescia \\ Via Valotti, 9, 25133 Brescia, Italy\\
%E-mail: paola.trebeschi@unibs.it
%}

\maketitle
\begin{abstract}
We study the free boundary problem for a plasma-vacuum interface in ideal incompressible magnetohydrodynamics. Unlike the classical statement when the vacuum magnetic field obeys the div-curl system of pre-Maxwell dynamics, to better understand the influence of the electric field in vacuum we do not neglect the displacement current in the vacuum region and consider the Maxwell equations for electric and magnetic fields. Under the necessary and sufficient stability condition for a planar interface found in \cite{T18}, we prove an energy a priori estimate for the linearized constant coefficient problem. The process of derivation of this estimate is based on various methods, including a secondary symmetrization of the vacuum Maxwell equations, the derivation of a hyperbolic evolutionary equation for the interface function and the construction of a degenerate Kreiss-type symmetrizer for an elliptic-hyperbolic problem for the total pressure.
\end{abstract}

\vspace{2mm}

\vspace{2mm}

\section{Introduction}
\label{s0}

%%%%%%%%%%%%%%%%%%%%%%%%%%%%%%%%%%%%%%%%%%%%%%%
%%%%%%%%%%%%%%%%%%%%%%%%%%%%%%%%%%%%%%%%%%%%%%%
%%%%%%%%%%%%%%%%%%%%%%%%%%%%%%%%%%%%%%%%%%%%%%%
%%%%%%%%%%%%%%%%%%%%%%%%%%%%%%%%%%%%%%%%%%%%%%%

 Plasma-vacuum interface problems are considered in the mathematical modeling of plasma confinement by magnetic fields in thermonuclear energy production (as in Tokamaks, Stellarators; see, e.g., \cite{BFKK}, \cite{Goed}). There are also important applications in astrophysics, where the plasma-vacuum interface problem can be used for modeling the motion of a star or the solar corona when magnetic fields are taken into account.

 Assume that the plasma-vacuum interface is described by $\Gamma (t)=\{F(t,x)=0\}$, and that
$\Omega^\pm(t)=\{F(t,x)\gtrless 0\}$ are the space-time domains occupied by the plasma and the vacuum respectively. Since $F$ is an unknown, this is a free-boundary problem.

 For the description of the motion in the plasma region let us introduce the unknowns $v=(v_1,v_2,v_3)$, $H=(H_1,H_2,H_3)$,  $q =p+\frac{1}{2}|{H} |^2$, $p$,  denoting respectively the velocity field , the magnetic field, the total pressure and the pressure. The unknowns in the vacuum region are $\mathcal{H}=(\mathcal{H}_1,\mathcal{H}_2,\mathcal{H}_3)$ and $\mathcal{E}=(\mathcal{E}_1,\mathcal{E}_2,\mathcal{E}_3)$ denoting respectively the magnetic field and the electric field.

 In the classical description of \cite{BFKK} (see also \cite{Goed}) the plasma is described by the  equations of ideal compressible magnetohydrodynamics (MHD)\footnote{Here we do not write out explicitly the compressible MHD equations because in the sequel we are going to consider the incompressible MHD equations in the plasma region.} whereas in the vacuum region one considers the so-called {\it pre-Maxwell dynamics}
\begin{equation}
\nabla \times \calH =0,\qquad {\rm div}\, {\calH}=0,\label{6}
\end{equation}
\begin{equation}
\nabla \times {\calE} =- \frac1{\varepsilon}\partial_t{\calH},\qquad {\rm div}\, \calE=0,\label{6'}
\end{equation}
where the positive constant $\varepsilon \ll 1$, being the ratio between a characteristic (average) speed of the plasma flow and the speed of light in vacuum, is a natural {\it small parameter} of the problem.
Notice that the equations \eqref{6}, \eqref{6'} are obtained from the Maxwell equations by neglecting the displacement current $(1/\varepsilon)\,\partial_t\calE$.
From \eqref{6'} the electric field $\calE$ is a secondary variable that may be  computed from the magnetic field ${ \calH}$; thus it is enough to consider \eqref{6} for the magnetic field $\calH$.

 The problem is completed by the boundary conditions at the free interface $\Gamma (t)$
\begin{subequations}\label{7}
\begin{align}
\frac{{\rm d}F }{{\rm d} t}=0,\quad
 [q]=0,\quad  H\cdot N=0,\label{8a}
 \\  \calH\cdot N=0 ,\label{8b}
\end{align}
 \end{subequations}
where $[q]=q_{|\Gamma}-\frac12|\calH_{|\Gamma}|^2$ denotes the jump of the total pressure across the interface and $N=\nabla F$. The
first condition in \eqref{8a} (where ${\rm d}/{\rm d}t= \partial_t +v\cdot\nabla$ denotes the material derivative) means that the interface moves with the velocity of plasma particles at the boundary.

 In \cite{T10} a basic energy a priori estimate in Sobolev spaces for the linearized plasma-vacuum interface problem was proved under the non-collinearity condition
\begin{equation}\label{non-coll}
H\times\calH\not=0 \qquad\mbox{on}\ \Gamma(t),
\end{equation}
satisfied on the interface for the unperturbed flow. Under the non-collinearity condition \eqref{non-coll} satisfied at the initial time, the well-posedness of the nonlinear problem was proved by Secchi and Trakhinin \cite{ST1,ST2} for compressible MHD equations in plasma region, and by Sun, Wang and Zhang \cite{SWZ2}  for incompressible MHD equations.

 The linearized stability of the relativistic case was first addressed by Trakhinin in \cite{T12}, in the case of plasma expansion in vacuum.  For the non-relativistic problem, the linearized stability was studied in the papers \cite{CDS,CDS2,MT} (see also \cite{M1,Tcpaa}) by considering a model where, in the vacuum region, instead of the pre-Maxwell dynamics the displacement current is taken into account and the complete system of Maxwell equations for the electric and the magnetic fields is considered. The introduction of this model aimed at investigating the influence of the electric field in vacuum on the well-posedness of the problem, since in the classical pre-Maxwell dynamics such an influence is hidden.

For the relativistic plasma-vacuum problem, Trakhinin had shown in \cite{T12} the possible ill-posedness in the presence of a {\it sufficiently strong} vacuum electric field. Since relativistic effects play a rather passive role in the analysis of  \cite{T12}, it was natural to expect a similar bahavior for the non-relativistic problem. In fact, it was shown in \cite{CDS,MT} that a {\it sufficiently weak} vacuum electric field precludes ill-posedness and gives the well-posedness of the linearized problem, thus somehow justifying the practice of neglecting the displacement current in the classical pre-Maxwell formulation when the vacuum electric field is weak enough. Such smallness hypothesis is not required for linear well-posedness in the two dimensional case, see \cite{CDS2}.

The results in \cite{CDS,MT,T12} induce a natural question: how strong the vacuum electric field has to be in order to enforce ill-posedness. The answer to this question has been given by Trakhinin in the recent paper \cite{T18} on the incompressible plasma-vacuum problem. In \cite{T18} the author analyses the linearized problem for the incompressible MHD equations in the plasma region and the Maxwell equations in the vacuum region, and obtains a necessary and sufficient condition for the violent instability of a planar plasma-vacuum interface (the opposite of this condition is given in \eqref{stabcond}).
In particular, it is shown that as the unperturbed plasma and vacuum magnetic fields are collinear (i.e. when \eqref{non-coll} is violated), any nonzero unperturbed vacuum electric field makes the planar interface violently unstable. This shows the necessity of the corresponding non-collinearity condition \eqref{non-coll} for well-posedness and a crucial role of the vacuum electric field in the evolution of a plasma-vacuum interface.

In the present paper we study the incompressible plasma-vacuum interface problem (the same problem as in \cite{T18}) and show that under the stability condition \eqref{stabcond} the linearized constant coefficient problem admits an energy a priori estimate, showing the stability of the planar plasma-vacuum interface.

%\textcolor{violet}{
%The proof of the a priori estimate requires various steps, including a secondary symmetrization of the vacuum Maxwell equations, the derivation of a hyperbolic evolutionary equation for the interface function and the construction of a degenerate Kreiss-type symmetrizer for an elliptic-hyperbolic problem  for the pressure.}

The rest of the paper is organized as follows. In Section \ref{s1} we first present the nonlinear plasma-vacuum interface problem. In Section \ref{s2} we introduce the constant solution \eqref{constsol} and derive the linearized problem \eqref{27'}--\eqref{29'} about \eqref{constsol}, which will be studied  in the paper. The main Theorem \ref{mainth1} is stated in Section \ref{s3} and its proof is given in the following sections. In Section \ref{s4}, we first introduce a secondary symmetrization of the vacuum Maxwell equations. Then, by  integrations by parts and a careful analysis of boundary integrals we prove an energy estimate showing that the solution in the interior domain may be controlled by the front function. In Section \ref{s5} we derive a hyperbolic evolution equation for the front function and obtain from it that the front function is estimated by the traces of the normal derivatives of the total pressures on both sides of the interface. In Section \ref{s6} we find an elliptic-hyperbolic problem for the total pressures, that we study in Section \ref{s8}, where we are able to get the crucial estimate of the traces of the normal derivatives of the total pressures, by constructing a degenerate Kreiss-type symmetrizer, and eventually close the desired estimate of the solution.

%%%%%%%%%%%%%%%%%%%%%%%%%%%%%%%%%%%%%%%%%%%%%%%
%%%%%%%%%%%%%%%%%%%%%%%%%%%%%%%%%%%%%%%%%%%%%%%
%%%%%%%%%%%%%%%%%%%%%%%%%%%%%%%%%%%%%%%%%%%%%%%
%%%%%%%%%%%%%%%%%%%%%%%%%%%%%%%%%%%%%%%%%%%%%%%

\section{The nonlinear problem}
\label{s1}

As in \cite{CDS,MT,Sec_sw,T18},  we assume that the free interface $\Gamma (t)$ has the form of a graph and the domains $\Omega^{\pm}(t)$ occupied by the plasma and the vacuum are unbounded:
\[
\Gamma (t)=\{F(t,x)=x_1-\varphi (t,x')\},\quad \Omega^\pm(t)=\{\pm (x_1- \varphi(t,x'))>0,\ x'\in \mathbb{R}^2\},\quad x'=(x_2,x_3).
\]
The plasma in the domain $\Omega^+(t)$ is assumed to be ideal and incompressible whereas in the vacuum region $\Omega^-(t)$
 we do not neglect the displacement current and consider the Maxwell equations. In a dimensionless form \cite{T18} the plasma-vacuum interface problem then reads:
\begin{equation}
\left\{
\begin{array}{ll}
{\rm div}\, v=0, & \\[6pt]
\displaystyle \frac{{\rm d}v}{{\rm d}t}-(H\cdot \nabla )H+\nabla q =0, & \\[12pt]
\displaystyle \frac{{\rm d}H}{{\rm d}t}-(H\cdot\nabla )v=0  &\quad \mbox{in}\ \Omega^+(t),
 \end{array}\right.
\label{13}
\end{equation}
\begin{equation}
\left\{
\begin{array}{l}
\varepsilon\,\partial_t\mathcal{H}+\nabla\times \mathcal{E} =0, \\[3pt]
\varepsilon\,\partial_t\mathcal{E}-\nabla\times \mathcal{H} =0\qquad  \mbox{in}\ \Omega^-(t),
 \end{array}\right.
\label{14}
\end{equation}
\begin{equation}
\left\{
\begin{array}{l}
\partial_t \varphi= v_N,\quad  q= \textstyle{\frac{1}{2}}\left(|\mathcal{H}|^2-|\mathcal{E}|^2\right), \\[3pt]
\mathcal{E}_{\tau_2}=\varepsilon\mathcal{H}_3\partial_t\varphi,\quad \mathcal{E}_{\tau_3}=-\varepsilon\mathcal{H}_2\partial_t\varphi \qquad \mbox{on}\ \Gamma (t),
 \end{array}\right.
\label{16}
\end{equation}
%\begin{equation}
%\label{17}
%\left\{
%\begin{array}{ll}
%{U} (0,{x})={U}_0({x}),\quad {x}\in \Omega^{+} (0),\quad
%\varphi (0,x')=\varphi_0(x'),\quad x'\in\mathbb{R}^2 , \\
%V(0,x)=
%V_0(x),\quad {x}\in \Omega^{-}(0),
%\end{array}
%\right.
%\end{equation}
where $U=(v,H)$, $V=(\mathcal{H},\mathcal{E})$. We denote
$v_N=v\cdot N$, $N= (1,-\partial_2\varphi ,-\partial_3\varphi )$ and $\mathcal{E}_{\tau_i}=\mathcal{E}_1\partial_i\varphi+\mathcal{E}_i$ ($i=2,3$). System \eqref{13}--\eqref{16} is supplemented with suitable initial conditions.
As for the case of compressible plasma flow in \cite{CDS,MT}, one can show that
\begin{equation}
{\rm div}\,{H}=0\quad \mbox{in}\ \Omega^+(t),\quad {\rm div}\,\mathcal{H}=0,\quad {\rm div}\,\mathcal{E}=0\quad \mbox{in}\ \Omega^-(t)
\label{18}
\end{equation}
and
\begin{equation}
{H}_N=0,\quad \mathcal{H}_N=0\quad \mbox{on}\ \Gamma (t),
\label{19}
\end{equation}
with ${H}_N=H\cdot N$ and $\mathcal{H}_N=\mathcal{H}\cdot N$. Equations \eqref{18}, \eqref{19} are the divergence and boundary constraints on the initial data, i.e., they hold for all $t>0$ if they are satisfied at $t=0$. The boundary conditions \eqref{16} are discussed for instance in \cite{MT,T18}.

\section{Linearized problem}

\label{s2}

Following \cite{T18}, we consider a solution $(U,V,\varphi )= (\widehat{U},\widehat{V},\sigma t)$ of  problem \eqref{13}--\eqref{16}
%(for $(U_0,V_0,\varphi_0 )=(\widehat{U},\widehat{V},0)$)
, where
\begin{equation}
\widehat{U}= (\hat{v},\widehat{H}),\quad \widehat{V}= (\widehat{\mathcal{H}},\widehat{\mathcal{E}}),
\label{constsol}
\end{equation}
with
\[
\begin{split}
& \hat{v}=(\sigma ,\hat{v}'),\quad \widehat{H}=(0,\widehat{H}'),\quad \widehat{\mathcal{H}}=(0,\widehat{\mathcal{H}}'),\quad
\widehat{\mathcal{E}}=(\widehat{\mathcal{E}}_1,\varepsilon\sigma\widehat{\mathcal{H}}_3,-\varepsilon\sigma
\widehat{\mathcal{H}}_2),
\\
& \hat{v}'=(\hat{v}_2,\hat{v}_3),\quad \widehat{H}'=(\widehat{H}_2,\widehat{H}_3),\quad \widehat{\mathcal{H}}'=(\widehat{\mathcal{H}}_2,\widehat{\mathcal{H}}_3),
\end{split}
\]
and $\hat{v}_k$, $\widehat{H}_k$, $\widehat{\mathcal{H}}_k$ ($k=2,3$), $\widehat{\mathcal{E}}_1$ and $\sigma$ are some constants. This solution describes a uniform flow with a planar interface moving with the constant velocity $\sigma$.

  Unlike the MHD system, the Maxwell equations {\it are not} Galilean invariant (they are Lorentz invariant), and we are not allowed to assume $\sigma =0$, as can be done without loss of generality in frequent situations when the equations on both sides of a planar interface are Galilean invariant.

Let us introduce the space-time domains
\[
\mathbb{R}^3_{\pm}=\{\pm x_1>0, x'\in\mathbb{R}^2\},\quad \Omega^{\pm}=\mathbb{R}\times \mathbb{R}^3_{\pm} =\{t\in\mathbb{R},\ x\in\mathbb{R}^3_{\pm}\}
\]
with the common boundary
\[
\omega = \mathbb{R}^3 = \{t\in\mathbb{R},\ x_1=0,\ x^\prime\in\mathbb{R}^2\}.
\]

We linearize problem \eqref{13}--\eqref{16} about the reference state \eqref{constsol} and introduce the change of independent variables
\[
\tilde{t}=t,\quad \tilde{x}_1= x_1-\sigma t,\quad \tilde{x}'=x'.
\]
We also denote
\[
 L=\partial_t +(\hat{v}'\cdot\nabla' ),\quad K=(\widehat{H}'\cdot\nabla' )\quad\mbox{and}\quad \nabla ' =(\partial_2,\partial_3).
\]
After dropping tildes for notational simplicity, the linearization gives the constant coefficient problem
\begin{equation}
\left\{
\begin{array}{ll}
{\rm div}\, v=f_0, & \\[3pt]
Lv-KH+\nabla q =f_1, & \\[3pt]
LH-Kv=f_2  &\quad \mbox{in}\ \Omega_+,
 \end{array}\right.
\label{26}
\end{equation}
\begin{equation}
\left\{
\begin{array}{l}
\varepsilon (\partial_t-\sigma\partial_1)\mathcal{H}+\nabla\times \mathcal{E} =f_3, \\[3pt]
\varepsilon (\partial_t-\sigma\partial_1)\mathcal{E}-\nabla\times \mathcal{H} =f_4\qquad  \mbox{in}\ \Omega_-,
 \end{array}\right.
\label{27}
\end{equation}
\begin{equation}
\left\{ \begin{array}{ll}
L\varphi= v_1+g_1, & \\[3pt]
q= \widehat{\mathcal{H}}_2(\mathcal{H}_2+\varepsilon \sigma \mathcal{E}_3) + \widehat{\mathcal{H}}_3(\mathcal{H}_3-\varepsilon \sigma \mathcal{E}_2) -\widehat{\mathcal{E}}_1\mathcal{E}_1 +g_2,  & \\[3pt]
\mathcal{E}_2=\varepsilon\widehat{\mathcal{H}}_3\partial_t\varphi-\widehat{\mathcal{E}}_1\partial_2\varphi +\varepsilon \sigma \mathcal{H}_3+g_3, & \\[3pt] \mathcal{E}_3=-\varepsilon\widehat{\mathcal{H}}_2\partial_t\varphi-\widehat{\mathcal{E}}_1\partial_3\varphi -\varepsilon \sigma \mathcal{H}_2+g_4 & \mbox{on}\ \omega,
\end{array}\right.
\label{28}
\end{equation}
for the perturbations $U=(v,H)$, $V=(\mathcal{H},\mathcal{E})$ and $\varphi$ (which are denoted by the same letters as the unknowns of the nonlinear problem and $q$ is the perturbation of the total pressure).

Following \cite{CDS}, we introduce the new unknowns
\begin{equation}
\begin{split}
\breve{\mathcal{H}} & =(\mathcal{H}_1,\mathcal{H}_2+\varepsilon\sigma \mathcal{E}_3,\mathcal{H}_3-\varepsilon\sigma \mathcal{E}_2),\\ \breve{\mathcal{E}} & = (\mathcal{E}_1,\mathcal{E}_2-\varepsilon\sigma \mathcal{H}_3,\mathcal{E}_3+\varepsilon\sigma\mathcal{H}_2).
\end{split}
\label{30}
\end{equation}
This is nothing else than the use of the nonrelativistic version of the Joules-Bernoulli equations (see, e.g., \cite{Sed}). In fact, in \cite{CDS} a more involved (``curved'') variant of \eqref{30} was applied  for showing that the corresponding plasma-vacuum interface problem for compressible MHD has a correct number of boundary conditions. In fact, the arguments in \cite{CDS} take also place for our case of incompressible MHD, and the number of boundary conditions in \eqref{28} is correct regardless of the sign of the interface speed $\sigma$.

After making the change of unknowns \eqref{30} and dropping breves the boundary conditions \eqref{28}  coincide with their form for $\sigma =0$. The Maxwell equations \eqref{27} can be written as
\begin{equation}\label{27^}
\varepsilon B_0\partial_t\breve{V} +\begin{pmatrix} \nabla\times \breve{\mathcal{E}}\\ -\nabla\times \breve{\mathcal{H}}\end{pmatrix} =0,
\end{equation}
where $\breve{V}=(\breve{\mathcal{H}},\breve{\mathcal{E}})$ and the matrix
\[
B_0=\frac{1}{1-\varepsilon^2\sigma^2}
\begin{pmatrix}
1-\varepsilon^2\sigma^2 & 0 & 0 & 0 & 0 & 0 \\
0 & 1 & 0 & 0 & 0 & -\varepsilon\sigma \\
0 & 0 & 1 & 0 & \varepsilon\sigma & 0 \\
0 & 0 & 0 & 1-\varepsilon^2\sigma^2 & 0 & 0\\
0 & 0& \varepsilon\sigma & 0 & 1 & 0\\
0 & -\varepsilon\sigma & 0 & 0 & 0& 1
\end{pmatrix}
\]
is found from the relation $V=B_0\breve{V}$. The matrix $B_0$ is symmetric and since $B_0|_{\varepsilon =0}=I_6$, in the {\it nonrelativistic limit} $\epsilon\to0$ we have $B_0>0$, i.e., system \eqref{27^} as well as the original system \eqref{27} is symmetric hyperbolic (here and below $I_k$ is the unit matrix of order $k$). In the nonrelativistic limit, $\varepsilon^2 \sigma^2 \rightarrow  0$ and spectral properties of the above constant coefficient hyperbolic system coincide for $\sigma =0$ and $\sigma \neq 0$ (for nonrelativistic speeds $\sigma$).
As in \cite{T18}, without loss of generality we may thus assume that $\sigma =0$. In fact, the below analysis of problem \eqref{26}--\eqref{28} with $\sigma =0$ just becomes a little bit more technically involved if instead of system \eqref{27} with $\sigma = 0$ we consider system \eqref{27^}, but there are no principal differences between the cases $\sigma \neq 0$ and $\sigma =0$.

Following \cite{CDS,MT,MTTpv}, we can reduce problem \eqref{26}--\eqref{28} to that with $f_0=0$, the homogeneous Maxwell equations ($f_3=f_4=0$), the homogeneous boundary conditions ($g_1=\ldots =g_4=0$) and the homogeneous divergence and boundary constraints (the linearizations of \eqref{18} and \eqref{19})
\begin{equation}
{\rm div}\,H=0\quad\mbox{in}\ \Omega^+,\quad {\rm div}\,{\mathcal{H}}=0,\quad {\rm div}\,\mathcal{E}=0\quad\mbox{in}\ \Omega^-,
\label{103}
\end{equation}
\begin{equation}
{H}_{1}=K\varphi ,\quad
{\mathcal{H}}_{1}=\mathcal{K}\varphi \quad\mbox{on}\ \omega,
\label{104}
\end{equation}
where $\mathcal{K} =(\mathcal{\widehat{H}}'\cdot\nabla' )$ (for zero initial data, \eqref{103} and \eqref{104} are automatically satisfied by the solutions of the reduced problem). To avoid overloading the paper we just refer the reader to \cite{CDS,MT,MTTpv} and do not describe here the process of partial homogenization of problem \eqref{26}--\eqref{28}.

For our subsequent analysis of problem \eqref{26}--\eqref{28} it will be also convenient to reflect the vacuum region $\Omega^-$ into the plasma domain $\Omega^+$, i.e., to make the change of variable $\check{x}_1=-x_1$ in \eqref{27}. Dropping checks, setting $\Omega:=\Omega^+$ and assuming that $\sigma =0$, $f_0=0$, $f_3=f_4=0$, $g_1=\ldots =g_4=0$, we obtain the following problem which is our main interest in this paper:
\begin{equation}
\left\{
\begin{array}{ll}
{\rm div}\, v=0, & \\[3pt]
Lv-KH+\nabla q =f_1, &  \\[3pt]
LH-Kv=f_2  & \\[3pt]
\varepsilon\, \partial_t\mathcal{H}+\nabla^-\times \mathcal{E} =0, \\[3pt]
\varepsilon \,\partial_t\mathcal{E}-\nabla^-\times \mathcal{H} =0\qquad \ \mbox{in}\ \Omega,
 \end{array}\right.
\label{27'}
\end{equation}
\begin{equation}
\left\{ \begin{array}{ll}
L\varphi= v_1, & \\[3pt]
q= \widehat{\mathcal{H}}_2\mathcal{H}_2+\widehat{\mathcal{H}}_3\mathcal{H}_3 -\widehat{\mathcal{E}}_1\mathcal{E}_1,  & \\[3pt]
\mathcal{E}_2=\varepsilon\widehat{\mathcal{H}}_3\partial_t\varphi-\widehat{\mathcal{E}}_1\partial_2\varphi , & \\[3pt] \mathcal{E}_3=-\varepsilon\widehat{\mathcal{H}}_2\partial_t\varphi-\widehat{\mathcal{E}}_1\partial_3\varphi  & \quad \mbox{on}\ \omega ,
\end{array}\right.
\label{28'}
\end{equation}
\begin{equation}
(U,V, \varphi)=0 \qquad \mbox{for}\ t<0,
\label{29'}
\end{equation}
where $\nabla^-=(-\partial_1,\partial_2,\partial_3)$ is the ``reflected'' operator of gradient.
We assume that the source terms $f_1$ and $f_2$ vanish in the past (i.e. for $t< 0$).

We remark that the solutions of problem \eqref{27'}--\eqref{29'} satisfy (cf. \eqref{103}, \eqref{104})
\begin{equation}
{\rm div}\,H=0,\quad {\rm div}^-\,{\mathcal{H}}=0,\quad {\rm div}^-\,\mathcal{E}=0\quad\mbox{in}\ \Omega,
\label{103'}
\end{equation}
\begin{equation}
{H}_{1}=K\varphi ,\quad
{\mathcal{H}}_{1}=\mathcal{K}\varphi \quad\mbox{on}\ \omega,
\label{104'}
\end{equation}
where the ``reflected'' divergence ${\rm div}^-\,a=-\partial_1a_1+\partial_2a_2+\partial_3a_3$ for any vector $a=(a_1,a_2,a_3)$.

\section{Main result}

\label{s3}

Before stating our main result we should introduce the weighted Sobolev spaces $H^m_{\gamma}(\Omega)$ and $H^m_{\gamma}(\omega)$, where $H^0_{\gamma}:=L^2_{\gamma}$,
$L^2_{\gamma}:=e^{\gamma t}L^2$, $H^m_{\gamma}:=e^{\gamma t}H^m$, with $\gamma \geq 1$, and the usual Sobolev spaces $H^m(\Omega)$ and $H^m(\omega)$ are equipped with the ({\it weighted}) norms
\begin{equation}\label{int_norm}
\nt u\nt^2_{m,\gamma:= \sum_{|\beta |\leq   m}\gamma^{2(m-|\beta |)} \| \partial^\beta  u\|^2_{L^2(\Omega)} \quad\mbox{and}\quad \|v\|^2_{m,\gamma} := \sum_{|\alpha|\leq m }
  \gamma^{2(m-|\alpha|)} \|\partial^\alpha_{\rm tan}  v \|^2_{L^2(\omega)}}
\end{equation}
respectively ($\partial_{\rm tan}^{\alpha}:=\partial_t^{\alpha_0}  \partial_{x_2}^{\alpha_2} \partial_{x_3}^{\alpha_3}$, with $\alpha =(\alpha_0, \alpha_2,\alpha_3)\in \mathbb{N}^3$). That is, the spaces $H^m_{\gamma}(\omega)$ and $H^m_{\gamma}(\Omega)$ are equipped with the norms
\[
\|u\|^2_{H^m_\gamma(\Omega ) :=\nt e^{-\gamma t}u\nt^2_{m,\gamma}\quad
\mbox{and}\quad\|v\|_{H^m_\gamma(\omega)} :=\|e^{-\gamma t}v\|_{m,\gamma}}
\]
for integer numbers $m$ and real $\gamma\ge 1$. Since in what follows we will also need to consider negative order Sobolev norms for functions on $\omega\equiv\mathbb R^3$, we recall that for any real order $m\in\mathbb R$ the Sobolev space $H^m(\omega)$ can be defined as the set of tempered distributions $v$ on $\omega$, making finite the (weighted) norm
\begin{equation}\label{frac_norm}
\Vert v\Vert^2_{m,\gamma}:=(2\pi)^{-2}\int_{\mathbb R^3}(\gamma^2+\vert\xi\vert^2)^m\vert\hat{v}(\xi)\vert^2d\xi
\end{equation}
being $\hat{v}=\hat{v}(\xi)$ the Fourier transform of $v$. In view of Plancherel's theorem, formula \eqref{frac_norm} is in agreement with the above definition of Sobolev norm with positive integer $m$, see the second formula in \eqref{int_norm}.

Observe that in terms of the weighted norms the trace estimate in $H^m$ reads
\begin{equation}\label{trace}
\|u_{|\omega}\|^2_{m,\gamma}\leq \frac{C}{\gamma}
\nt u\nt^2_{m+1,\gamma}.
\end{equation}

We are now in a position to state the main result of the paper.

\begin{theorem}\label{mainth1}
For every given planar plasma-vacuum interface described by the constant solution \eqref{constsol} and satisfying the stability condition
\begin{equation}\label{stabcond}
\widehat{\mathcal{E}}_1^{\,2}< \frac{|\widehat{H}|^2 +|\widehat{\mathcal{H}}|^2-\sqrt{\big(|\widehat{H}|^2 +|\widehat{\mathcal{H}}|^2\big)^2-4|\widehat{H}\times\widehat{\mathcal{H}}|^2}}{2},
\end{equation}
there exist constants $\varepsilon_*>0$ and $C>0$ such that for all $0<\varepsilon <\varepsilon_*$, $\gamma \geq 1$, any solution $((U,V),\varphi)\in L^2_\gamma(\Omega)\times H^1_\gamma(\omega)$ of problem \eqref{27'}--\eqref{29'}, with source term $f=(f_1,f_2)\in H^3_\gamma(\Omega)$ vanishing in the past, obeys the a priori estimate
\begin{equation}\label{est0}
\|(U,V)\|^2_{L^2_\gamma(\Omega)}+\| (U,\mathcal{H}_1,\mathcal{E}_2,
\mathcal{E}_3)|_{\omega}\|^2_{L^2_\gamma (\omega)} +\| \varphi\|^2_{H^1_\gamma(\omega)} \leq
 \frac{C}{\gamma^4} \|f\|^2_{H^3_\gamma(\Omega)}\,.
\end{equation}
\end{theorem}

\begin{remark}
{\rm
In the above theorem the assertion about the existence of a (small) value $\varepsilon_*$ just means that the necessary and sufficient neutral stability condition \eqref{stabcond} found in  \cite{T18} is  valid in the {\it nonrelativistic limit} $\varepsilon \rightarrow 0$.
}
\label{r1}
\end{remark}

It will be more convenient to prove Theorem \ref{mainth1} after its reformulation in terms of the exponentially weighted unknowns
\begin{equation}\label{change}
\bar{U}:=e^{-\gamma t}U,\quad \bar{V}:=e^{-\gamma t}V,\quad \bar{q}:=e^{-\gamma t}q,\quad\bar{\varphi}:=e^{-\gamma t}\varphi.
\end{equation}
We first restate problem \eqref{27'}--\eqref{29'} in terms the unknowns \eqref{change}:
\begin{equation}
\left\{
\begin{array}{ll}
{\rm div}\, \bar{v}=0, & \\[3pt]
L_{\gamma}\bar{v}-K\bar{H}+\nabla \bar{q} =\bar{f}_1, &  \\[3pt]
L_{\gamma}\bar{H}-K \bar{v}=\bar{f}_2  & \\[3pt]
\varepsilon(\gamma I +\partial_t)\bar{\mathcal{H}}+\nabla^-\times \bar{\mathcal{E}} =0, \\[3pt]
\varepsilon (\gamma I +\partial_t)\bar{\mathcal{E}}-\nabla^-\times \bar{\mathcal{H}} =0\qquad \ \mbox{in}\ \Omega,
 \end{array}\right.
\label{27''}
\end{equation}
\begin{equation}
\left\{ \begin{array}{ll}
L_{\gamma}\bar{\varphi}=  \bar{v}_1, & \\[3pt]
\bar{q}= \widehat{\mathcal{H}}_2\bar{\mathcal{H}}_2+\widehat{\mathcal{H}}_3\bar{\mathcal{H}}_3 -\widehat{\mathcal{E}}_1\bar{\mathcal{E}}_1,  & \\[3pt]
\bar{\mathcal{E}}_2=\varepsilon\widehat{\mathcal{H}}_3(\gamma I +\partial_t)\bar{\varphi}-\widehat{\mathcal{E}}_1\partial_2\bar{\varphi} , & \\[3pt] \bar{\mathcal{E}}_3=-\varepsilon\widehat{\mathcal{H}}_2(\gamma I +\partial_t)\bar{\varphi}-\widehat{\mathcal{E}}_1\partial_3\bar{\varphi}  & \quad \mbox{on}\ \omega ,
\end{array}\right.
\label{28''}
\end{equation}
\begin{equation}
(\bar{U},\bar{V}, \bar\varphi)=0 \qquad \mbox{for}\ t<0,
\label{29''}
\end{equation}
where $L_{\gamma}=\gamma I  +L$, with the identity operator $I$; $\bar{f}_i=e^{-\gamma t}f_i$, $i=1,2$, $\bar{v}=e^{-\gamma t}v$, etc.  For the new unknowns \eqref{change} equations \eqref{103'} and \eqref{104'} remain unchanged. Theorem \ref{mainth1} then admits the following equivalent formulation.

\begin{theorem}\label{mainth2}
For every given planar plasma-vacuum interface described by the constant solution \eqref{constsol} and satisfying the stability condition \eqref{stabcond}
there exist constants $\varepsilon_*>0$ and $C>0$ such that for all $0<\varepsilon <\varepsilon_*$, $\gamma \geq 1$, any solution $((\bar{U},\bar{V}),\bar{\varphi})\in L^2(\Omega)\times H^1(\omega)$ of problem \eqref{27''}--\eqref{29''}, with source term $\bar{f}=(\bar{f}_1,\bar{f}_2)\in H^3(\Omega)$ vanishing in the past, obeys the a priori estimate
\begin{equation}\label{est1.1}
\|(\bar{U},\bar{V})\|^2_{L^2(\Omega)}+\| (\bar{U},\bar{\mathcal{H}}_1,\bar{\mathcal{E}}_2,\bar{\mathcal{E}}_3)|_{\omega}\|^2_{L^2(\omega)} +\| \bar{\varphi}\|^2_{1,\gamma} \leq
 \frac{C}{\gamma^4}\nt\bar{f}\nt^2_{3,\gamma}\,.
\end{equation}
\end{theorem}

\section{Estimate of the interior unknowns through the interface \\ function}
\label{s4}

In order to simplify the notations, from now on we drop bars in problem \eqref{27''}--\eqref{29''} and the desired estimate \eqref{est1.1}. We first rewrite system \eqref{27''} as follows:
\begin{align}
& {\rm div}\, v=0, \label{v} \\
& \gamma U+\partial_tU +{A}_2\partial_2U+{A}_3\partial_3U +\begin{pmatrix}\nabla q \\ 0\end{pmatrix}=f,  \label{U} \\
&  \gamma V+\partial_tV +\varepsilon^{-1}\sum_{j=1}^{3}B_j\partial_jV =0  \qquad\qquad\ \mbox{in}\ \Omega, \label{V}
\end{align}
where
\[
A_k=\left( \begin{array}{cc} \hat{v}_k & -\widehat{H}_k\\
-\widehat{H}_k & \hat{v}_k\\
\end{array} \right)\otimes I_3,\quad k=1,2,3,\quad
B_1=\left(\begin{array}{cccccc}
0 & 0 & 0& 0 & 0 & 0 \\
0 & 0 & 0& 0 & 0 & 1 \\
0 & 0 & 0& 0 & -1 & 0 \\
0 & 0 & 0& 0 & 0 & 0 \\
0 & 0 & -1& 0 & 0 & 0 \\
0 & 1 & 0& 0 & 0 & 0
\end{array} \right),
\]
\[
B_2=\left(\begin{array}{cccccc}
0 & 0 & 0& 0 & 0 & 1 \\
0 & 0 & 0& 0 & 0 & 0 \\
0 & 0 & 0& -1 & 0 & 0 \\
0 & 0 & -1& 0 & 0 & 0 \\
0 & 0 & 0& 0 & 0 & 0 \\
1 & 0 & 0& 0 & 0 & 0
\end{array} \right),\quad
B_3=\left(\begin{array}{cccccc}
0 & 0 & 0& 0 & -1 & 0 \\
0 & 0 & 0& 1 & 0 & 0 \\
0 & 0 & 0& 0 & 0 & 0 \\
0 & 1 & 0& 0 & 0 & 0 \\
-1 & 0 & 0& 0 & 0 & 0 \\
0 & 0 & 0& 0 & 0 & 0
\end{array} \right).
\]

The crucial role in the argument below of the energy method will be played by the so-called secondary symmetrization of the vacuum Maxwell equations proposed in \cite{T12}. Following \cite{T12} (see also \cite{CDS,MT,ST1})  and using the last two divergences in \eqref{103'}, we equivalently rewrite system \eqref{V} as the symmetric system
\begin{equation}
\mathcal{B}_0(\gamma I+\partial_t)V +\varepsilon^{-1}\sum_{j=1}^{3}\mathcal{B}_j\partial_jV =0  \quad\ \mbox{in}\ \Omega, \label{V'}
\end{equation}
where
\[
\mathcal{B}_0=\left(\begin{array}{cccccc}
1 & 0 & 0& 0 & \nu_3 & -\nu_2 \\
0 & 1 & 0& -\nu_3 & 0 & \nu_1 \\
0 & 0 & 1& \nu_2 & -\nu_1 & 0 \\
0 & -\nu_3 & \nu_2& 1 & 0 & 0 \\
\nu_3 & 0 & -\nu_1& 0 & 1 & 0 \\
-\nu_2 & \nu_1 & 0& 0 & 0 & 1
\end{array} \right),\quad
\mathcal{B}_1=
\left(\begin{array}{cccccc}
-\nu_1 & -\nu_2 & -\nu_3& 0 & 0 & 0 \\
-\nu_2 & \nu_1 & 0& 0 & 0 & 1 \\
-\nu_3 & 0 & \nu_1& 0 & -1 & 0 \\
0 & 0 & 0& -\nu_1 & -\nu_2 & -\nu_3 \\
0 & 0 & -1& -\nu_2 & \nu_1 & 0 \\
0 & 1 & 0& -\nu_3 & 0 & \nu_1
\end{array} \right),
\]
\[
\mathcal{B}_2=
\left(\begin{array}{cccccc}
-\nu_2 & \nu_1 & 0& 0 & 0 & 1 \\
\nu_1 & \nu_2 & \nu_3& 0 & 0 & 0 \\
0 & \nu_3 & -\nu_2& -1 & 0 & 0 \\
0 & 0 & -1& -\nu_2 & \nu_1 & 0 \\
0 & 0 & 0& \nu_1 & \nu_2 & \nu_3 \\
1 & 0 & 0& 0 & \nu_3 & -\nu_2
\end{array} \right),\quad
\mathcal{B}_3=
\left(\begin{array}{cccccc}
-\nu_3 & 0 & \nu_1& 0 & -1 & 0 \\
0 & -\nu_3 & \nu_2& 1 & 0 & 0 \\
\nu_1 & \nu_2 & \nu_3& 0 & 0 & 0 \\
0 & 1 & 0& -\nu_3 & 0 & \nu_1 \\
-1 & 0 & 0& 0 & -\nu_3 & \nu_2 \\
0 & 0 & 0& \nu_1 & \nu_2 & \nu_3
\end{array} \right),
\]
and  $\nu_i$ ($i=1,2,3$) are arbitrary constants satisfying the hyperbolicity condition $\mathcal{B}_0>0$, i.e.,
\begin{equation}
|\nu |<1,\label{21'}
\end{equation}
with $\nu = (\nu_1, \nu_2, \nu_3)$. Because of the reflection of the vacuum region made above the matrix $\mathcal{B}_1$ has here the opposite sign in comparison with that in \cite{CDS,ST1}.

We now make the same choice of the constant vector $\nu$ as in \cite{MT}:
\[
\nu =\varepsilon \hat{v}=\varepsilon (0,\hat{v}_2,\hat{v}_3).
\]
For this choice the hyperbolicity condition \eqref{21'} holds in the nonrelativistic limit. Using standard arguments of the energy method and taking into account $\mathcal{B}_0>0$ and the incompressibility condition \eqref{v}, for systems \eqref{U} and \eqref{V'} we deduce the energy inequality
\begin{equation}\label{ei}
\gamma\|U\|^2_{L^2(\Omega)}+
\gamma\mathcal{I}
+\int\limits_{\omega}\mathcal{Q}\,{\rm d}x'{\rm d}t\leq \frac{C}{\gamma}\|{f}\|^2_{L^2(\Omega)},
\end{equation}
where
\[
\mathcal{I}=\int\limits_{\Omega}(\mathcal{B}_0V\cdot V){\rm d}x{\rm d}t
\]
and
\[
\begin{split}
\mathcal{Q}  &=-qv_1|_{\omega} -\frac{1}{2\varepsilon}(\mathcal{B}_1V\cdot V)|_{\omega} \\ & =\left.\left\{ -qv_1 +\mathcal{H}_1(\hat{v}_2\mathcal{H}_2+\hat{v}_3\mathcal{H}_3)+\mathcal{E}_1(\hat{v}_2\mathcal{E}_2+\hat{v}_3\mathcal{E}_3)+
\varepsilon^{-1}(\mathcal{H}_3\mathcal{E}_2-\mathcal{H}_2\mathcal{E}_3)\right\}\right|_{\omega}.
\end{split}
\]
Here and below $C >0$ is a constant independent of $\gamma$.

Using the boundary conditions \eqref{28''} and the second condition in \eqref{104'}, after  some algebra we get
\[
\mathcal{Q}=\hat{\mu}\left.\left\{ \mathcal{E}_1(\gamma\varphi +\partial_t\varphi ) +\varepsilon^{-1}(\mathcal{H}_2\partial_3\varphi -\mathcal{H}_3\partial_2\varphi ) \right\}\right|_{\omega},
\]
where $\hat{\mu}=\widehat{\mathcal E}_1+\varepsilon\hat{v}_2\widehat{\mathcal{H}}_3-\varepsilon\hat{v}_3\widehat{\mathcal{H}}_2$. Using then the fourth equation in \eqref{V} restricted to the boundary, we rewrite $\mathcal{Q}$ as follows:
\[
\mathcal{Q}=2\gamma\hat{\mu}\varphi\mathcal{E}_1|_{\omega}+ \partial_t\left(\hat{\mu}\varphi \mathcal{E}_1|_{\omega} \right)-\partial_2\left(\varepsilon^{-1}\hat{\mu}\varphi\mathcal{H}_3|_{\omega} \right)  + \partial_3\left(\varepsilon^{-1}\hat{\mu}\varphi \mathcal{H}_2|_{\omega}\right) .
\]
The substitution of the last formula into \eqref{ei} gives, cf. \cite{CDS,MT},
\begin{equation}\label{ei'}
\gamma\|U\|^2_{L^2(\Omega)}+
\gamma\mathcal{I}
\leq -2\gamma\int\limits_{\omega}\hat{\mu}\varphi\mathcal{E}_1|_{\omega}{\rm d}x'{\rm d}t+\frac{C}{\gamma}\|{f}\|^2_{L^2(\Omega)}.
\end{equation}

From the boundary integral in \eqref{ei'} we pass to the volume integral, use the third divergence in \eqref{103'}, integrate by parts and apply the Young inequality with an arbitrary positive constant $\delta$:
\begin{equation}\label{ei''}
\begin{split}
-2\gamma\int\limits_{\omega}\hat{\mu}\varphi \mathcal{E}_1|_{\omega}{\rm d}x'{\rm d}t & =2\gamma\int\limits_{\Omega} \hat{\mu}\chi \varphi \,\partial_1\mathcal{E}_1{\rm d}x{\rm d}t =
+2\gamma\int\limits_{\Omega}\hat{\mu}\chi \varphi \left(\partial_2\mathcal{E}_2+\partial_3\mathcal{E}_3\right){\rm d}x{\rm d}t\\
& =-2 \gamma\int\limits_{\Omega}\hat{\mu}\chi \left(\mathcal{E}_2\partial_2\varphi +\mathcal{E}_3\partial_3\varphi\right){\rm d}x{\rm d}t
 \leq \gamma C\Big(\delta \mathcal{I} +\frac{1}{\delta}\| \varphi\|^2_{1,\gamma}\Big),
\end{split}
\end{equation}
where the lifting function $\chi (x_1)\in C_0^{\infty}(\mathbb{R}_+)$ can be taken, for example, such that $\chi =1$ on $[0,1/2]$ and $\chi =0$ on $[1,\infty)$. Choosing $\delta$ small enough and taking into account that $\mathcal{B}_0>0$, from \eqref{ei'} and \eqref{ei''} we derive the following estimate of the interior unknowns $U$ and $V$
through the interface function $\varphi$ (and the source term $f$) announced in the title of this section:
\begin{equation}\label{UV}
\gamma\|(U,V)\|^2_{L^2(\Omega)}\leq C\Big(\gamma\| \varphi\|^2_{1,\gamma} + \frac{1}{\gamma}\|{f}\|^2_{L^2(\Omega)}\Big).
\end{equation}

\section{Hyperbolic evolution equation and estimate for the interface function}

\label{s5}

Following ideas of \cite{SWZ1,SWZ2}, we now deduce an evolution equation for the interface function. From the first condition in \eqref{104'}, the second scalar equation in \eqref{27''} and the first boundary condition in \eqref{28''} we obtain that
\begin{equation}\label{phi1}
(L_{\gamma}^2-K^2)\varphi +\partial_1 q = f_{1,1}\quad \mbox{on}\ \omega ,
\end{equation}
where $f_{1,1}$ is the first component of the source term $f_1=(f_{1,1},f_{1,2},f_{1,3})$. Let us introduce the perturbation of the vacuum total pressure whose trace appears in the right-hand side of the second boundary condition in \eqref{28''}:
\[
q^-:=\widehat{\mathcal{H}}_2\mathcal{H}_2+\widehat{\mathcal{H}}_3\mathcal{H}_3 -\widehat{\mathcal{E}}_1\mathcal{E}_1.
\]
Then, we rewrite \eqref{phi1} as
\begin{equation}\label{phi2}
(L_{\gamma}^2-K^2)\varphi - \widehat{\mathcal{E}}_1\partial_1\mathcal{E}_1
+\widehat{\mathcal{H}}_2\partial_1\mathcal{H}_2+\widehat{\mathcal{H}}_3\partial_1\mathcal{H}_3+\partial_1 q -\partial_1 q^-
= f_{1,1}\quad \mbox{on}\ \omega .
\end{equation}

Using the third divergence in \eqref{103'} and the last two boundary conditions in \eqref{28''},  we get
\begin{equation}\label{phi3}
\partial_1\mathcal{E}_1|_{\omega}=- \widehat{\mathcal{E}}_1\Delta'\varphi +\varepsilon (\gamma I +\partial_t)\mathcal{K}^{\bot}{\varphi},
\end{equation}
where $\Delta'=\partial_2^2+\partial_3^2$ is the tangential Laplacian and $\mathcal{K}^{\bot}=\widehat{\mathcal{H}}_3\partial_2-\widehat{\mathcal{H}}_2\partial_3$. From the last two equations in \eqref{27''} we have
\begin{equation}\label{phi4}
\partial_1\mathcal{H}_2=- \varepsilon (\gamma I +\partial_t)\mathcal{E}_3-\partial_2\mathcal{H}_1\quad\mbox{and}\quad
\partial_1\mathcal{H}_3= \varepsilon (\gamma I +\partial_t)\mathcal{E}_2+\partial_3\mathcal{H}_1.
\end{equation}
Restricting \eqref{phi4} to $\omega$ and using the last boundary conditions in \eqref{28''} as well as the second condition in \eqref{104'}, we obtain
\begin{equation}\label{phi5}
\begin{split}
(\widehat{\mathcal{H}}_2\partial_1\mathcal{H}_2+\widehat{\mathcal{H}}_3\partial_1\mathcal{H}_3)|_{\omega} & = -\mathcal{K}^2\varphi+\varepsilon (\gamma I +\partial_t)
(\mathcal{H}_3\mathcal{E}_2-\mathcal{H}_2\mathcal{E}_3)|_{\omega} \\
& = -\mathcal{K}^2\varphi +\varepsilon^2 |\widehat{\mathcal{H}}'|^2(\gamma I +\partial_t)^2\varphi-\varepsilon \widehat{\mathcal{E}}_1(\gamma I +\partial_t)\mathcal{K}^{\bot}{\varphi}.
\end{split}
\end{equation}

Substituting \eqref{phi3} and \eqref{phi5} into \eqref{phi2}, we get the desired evolution equation
\begin{equation}\label{eveq}
\mathcal{L}_{\gamma}\varphi = F\quad \mbox{on}\ \omega ,
\end{equation}
where
\[
\begin{split}
\mathcal{L}_{\gamma} & =L_{\gamma}^2-K^2-\mathcal{K}^2+\widehat{\mathcal{E}}_1^{\,2}\Delta' -2\varepsilon \widehat{\mathcal{E}}_1(\gamma I +\partial_t)\mathcal{K}^{\bot} +\varepsilon^2 |\widehat{\mathcal{H}}'|^2(\gamma I +\partial_t)^2,\\ F & =(\partial_1 q^- -\partial_1 q +f_{1,1})|_{\omega}.
\end{split}
\]
In the nonrelativistic setting $\varepsilon \ll 1$, the operator $\mathcal{L}_\gamma$ is {\it hyperbolic} provided that the operator $L_{\gamma}^2-K^2-\mathcal{K}^2+\widehat{\mathcal{E}}_1^{\,2}\Delta'$ does. At the same time, the principal part of the last operator is the operator
\[
P=L^2-K^2-\mathcal{K}^2+\widehat{\mathcal{E}}_1^{\,2}\Delta'.
\]
Considered that $  L=\partial_t +(\hat{v}'\cdot\nabla' )$ is a transport operator, the operator $P$ is hyperbolic if and only if the quadratic form
\[
Q(x,y)= \big(\widehat{H}_2x+\widehat{H}_3y\big)^2 +\big(\widehat{\mathcal{H}}_2x+\widehat{\mathcal{H}}_3y\big)^2-\widehat{\mathcal{E}}_1^{\,2}(x^2+y^2)>0.
\]
One can check that the latter is true if and only if the stability condition \eqref{stabcond} holds,  see \cite{T18}.

Let \eqref{stabcond} be fulfilled. Then, considering for a moment $F$ as a given right-hand side in \eqref{eveq}, we have a hyperbolic equation  for $\varphi$. As for the wave equation (written in terms of an ``exponentially weighted'' unknown $e^{-\gamma t}u$), from \eqref{eveq} we get the a priori estimate
\[
\gamma \|\varphi\|^2_{1,\gamma}\leq \frac{C}{\gamma}\|F\|^2_{L^2(\omega )}.
\]
Using for $f_{1,1}|_{\omega}$ the trace theorem, we come to the estimate
\begin{equation}\label{phi6}
\gamma \|\varphi\|^2_{1,\gamma}\leq C \left( \frac{1}{\gamma^2}\nt f\nt^2_{1,\gamma}+ \frac{1}{\gamma}\|\partial_1q|_{\omega}\|^2_{L^2(\omega )}
+ \frac{1}{\gamma}\|\partial_1q^-|_{\omega}\|^2_{L^2(\omega )} \right).
\end{equation}
In view of \eqref{phi6}, for ``closing'' estimate \eqref{UV} it remains to estimate the traces $\partial_1q|_{\omega}$ and $\partial_1q^-_{|\omega}$ through the source $f$.

\section{Elliptic-hyperbolic problem for the total pressures}
\label{s6}

Clearly, using the first divergence in \eqref{103'}, from \eqref{v} and \eqref{U} we derive the Poisson equation
\begin{equation}\label{q}
\Delta q={\rm div}\, f_1\quad \mbox{in}\ \Omega .
\end{equation}
From system \eqref{V} we derive the wave equation for each component of $V$. Hence, we obtain the wave equation for $q^-$:
\begin{equation}\label{q-}
\varepsilon^2(\gamma I +\partial_t)^2q^- -\Delta q^-=0 \quad \mbox{in}\ \Omega .
\end{equation}

It follows from \eqref{phi3} and \eqref{phi5} that
\begin{equation}\label{q-b}
\partial_1q^-=(\widehat{\mathcal{E}}_1^{\,2}\Delta' -\mathcal{K}^2)\varphi -2\varepsilon \widehat{\mathcal{E}}_1(\gamma I +\partial_t)\mathcal{K}^{\bot} \varphi +\varepsilon^2 |\widehat{\mathcal{H}}'|^2(\gamma I +\partial_t)^2\varphi\quad \mbox{on}\ \omega .
\end{equation}
By adding \eqref{phi1} and \eqref{q-b}, we get
\begin{equation}\label{qq+}
\partial_1q+\partial_1q^-=\mathcal{P}_{\gamma}\varphi+f_{1,1}\quad \mbox{on}\ \omega ,
\end{equation}
with
\[
\begin{split}
\mathcal{P}_{\gamma} & =K^2-\mathcal{K}^2-L_{\gamma}^2+\widehat{\mathcal{E}}_1^{\,2}\Delta' -2\varepsilon \widehat{\mathcal{E}}_1(\gamma I +\partial_t)\mathcal{K}^{\bot} +\varepsilon^2 |\widehat{\mathcal{H}}'|^2(\gamma I +\partial_t)^2 \\ & =\mathcal{L}_\gamma -2L_{\gamma}^2+2K^2.
\end{split}
\]
At the same time, we rewrite \eqref{eveq} as
\begin{equation}\label{qq-}
\partial_1q-\partial_1q^-=-\mathcal{L}_{\gamma}\varphi+f_{1,1}\quad \mbox{on}\ \omega .
\end{equation}

Applying $\mathcal{L}_\gamma$ to \eqref{qq+} and $\mathcal{P}_\gamma$ to \eqref{qq-} and adding the results, we exclude $\varphi$ from \eqref{qq+} and \eqref{qq-}:
\begin{equation}\label{qq}
\Sigma^-_{\gamma}\partial_1q +\Sigma^+_{\gamma}\partial_1q^-=\Sigma^-_{\gamma}f_{1,1} \quad \mbox{on}\ \omega ,
\end{equation}
where
\[
\begin{split}
&\Sigma^-_{\gamma}=\frac{\mathcal L_\gamma+\mathcal P_\gamma}{2}=\mathcal{L}_\gamma -L_{\gamma}^2+K^2=-\mathcal{K}^2+\widehat{\mathcal{E}}_1^{\,2}\Delta' -2\varepsilon \widehat{\mathcal{E}}_1(\gamma I +\partial_t)\mathcal{K}^{\bot} +\varepsilon^2 |\widehat{\mathcal{H}}'|^2(\gamma I +\partial_t)^2,\\
&\Sigma^+_{\gamma}=\frac{\mathcal L_\gamma-\mathcal P_\gamma}{2}=L_{\gamma}^2-K^2.
\end{split}
\]
Collecting \eqref{q}, \eqref{q-}, \eqref{qq} and the second boundary condition in \eqref{28''}, we get the following elliptic-hyperbolic problem for the total pressures $q$ and $q^-$:
\begin{equation}
\left\{
\begin{array}{l}
\Delta q={\rm div}\, f_1, \\
\varepsilon^2(\gamma I +\partial_t)^2q^- -\Delta q^-=0 \quad \mbox{in}\ \Omega,
\end{array}
\right.\label{01}
\end{equation}
\begin{equation}
\left\{
\begin{array}{l}
q-q^-=0,\\
\Sigma^-_{\gamma}\partial_1q +\Sigma^+_{\gamma}\partial_1q^-=\Sigma^-_{\gamma}f_{1,1} \quad \mbox{on}\ \omega .
\end{array}
\right. \label{02}
\end{equation}

In the subsequent analysis it will be more convenient to have fully homogeneous interior equations in \eqref{01}, i.e., the Laplace equation instead of the Poisson equation. Following \cite{MTT1}, we introduce the ``shift'' $\breve{q}$ satisfying the elliptic problem
\begin{align}
 \Delta \breve{q}={\rm div}\, f_1 \quad \mbox{in}\ \Omega, \label{03}\\
\breve{q} =\partial_1 \breve{q} +f_{1,1}\quad \mbox{on}\ \omega .\label{04}
\end{align}
Multiplying \eqref{03} by $\breve{q}$, integrating the result over $\Omega$ and using the boundary condition \eqref{04}, we get by standard arguments the estimate
\[
\|\partial_1\breve{q}|_{\omega}\|^2_{L^2(\omega )}+\frac{1}{2}\|\nabla \breve{q}\|^2_{L^2(\Omega )} \leq \|f_{1,1}|_{\omega}\|^2_{L^2(\omega )}+\frac{1}{2}\| f_1\|^2_{L^2(\Omega )}.
\]
Using again \eqref{04} and the trace theorem, we finally obtain the estimate
\begin{equation}\label{05}
\|\breve{q}|_{\omega}\|^2_{L^2(\omega )} +
\|\partial_1\breve{q}|_{\omega}\|^2_{L^2(\omega )} \leq C\nt f_1\nt^2_{1,\gamma}.
\end{equation}
Clearly, the tangential differentiation of problem \eqref{03}, \eqref{04} gives us also the estimate
\begin{equation}\label{06}
\|\breve{q}|_{\omega}\|^2_{m,\gamma} +
\|\partial_1\breve{q}|_{\omega}\|^2_{m,\gamma} \leq C\nt f_1\nt^2_{m+1,\gamma}, \quad \forall m\in\mathbb{N}.
\end{equation}
It follows from \eqref{05}, \eqref{06} and the elementary inequality
\[
\nt u\nt_{s,\gamma}\leq \frac{1}{\gamma^{r-s}}
\nt u\nt_{r,\gamma} \quad \mbox{for}\ r>s
\]
that
\begin{equation}\label{07}
 \|\nabla\breve{q}|_{\omega}\|^2_{L^2(\omega )} \leq C\nt f_1\nt^2_{2,\gamma}.
\end{equation}

We now introduce the ``shifted'' total pressure $q^+=q-\breve{q}$. It follows from \eqref{01}--\eqref{04} that $q^+$ and $q^-$ satisfy the problem
\begin{equation}
\left\{
\begin{array}{l}
\Delta q^+=0, \\
\varepsilon^2(\gamma I +\partial_t)^2q^- -\Delta q^-=0 \quad \mbox{in}\ \Omega,
\end{array}
\right.\label{01'}
\end{equation}
\begin{equation}
\left\{
\begin{array}{l}
q^+-q^-=\mathfrak{g}_2,\\
\Sigma^-_{\gamma}\partial_1q^+ +\Sigma^+_{\gamma}\partial_1q^-=\mathfrak{g}_1 \quad \mbox{on}\ \omega ,
\end{array}
\right. \label{02'}
\end{equation}
where
\begin{equation}\label{gg}
\mathfrak{g}_2=-\breve{q}|_{\omega},\quad \mathfrak{g}_1=\Sigma^-_{\gamma}f_{1,1}|_{\omega}-\Sigma^-_{\gamma}\partial_1\breve{q}|_{\omega}.
\end{equation}

\section{Construction of a symbolic symmetrizer for problem \eqref{01'}, \eqref{02'}}
\label{s8}

\subsection{A boundary value problem for the Fourier transforms}

%Our present goal is deriving an a priori estimate for the traces of $\nabla q$, $\nabla q^-$ and the front function $\varphi$.
We first apply a Fourier transform to problem \eqref{01'}, \eqref{02'} with respect to $x'=(x_2,x_3)$ and $t$, with the Fourier dual variables $\eta'=(\eta_2,\eta_3)$ and $\delta$ respectively. Let us also set
\[
\tau =\gamma +i\delta
\]
and
\[
\Lambda(\tau,\eta^\prime):=\sqrt{\vert\tau\vert^2+\eta^2}\quad\mbox{and}\quad\eta:=|\eta'|\,.
\]
For the Fourier transformed pressures
\begin{align*}
\tilde{q}^+(\delta ,x_1,\eta'):=\int\limits_{\mathbb R^3}e^{-i\delta t-i\eta^\prime\cdot x^\prime}q^+(t,x_1,x^\prime)dtdx^\prime\,,\\
\tilde{q}^-(\delta,x_1,\eta'):=\int\limits_{\mathbb R^3}e^{-i\delta t-i\eta^\prime\cdot x^\prime}q^-(t,x_1,x^\prime)dtdx^\prime
\end{align*}
from \eqref{01'}, \eqref{02'} we obtain the following problem
\begin{align}
\left\{
\begin{array}{l}
\displaystyle \frac{d^2\tilde{q}^+}{dx_1^{\,2}}-\eta^2\tilde{q}^+=0, \\[12pt]
\displaystyle \frac{d^2\tilde{q}^-}{dx_1^{\,2}}-(\eta^2 +\varepsilon^2\tau^2 )\tilde{q}^-=0,\quad x_1>0,
\end{array}
\right.\label{9}\\[6pt]
\left\{
\begin{array}{l}
\displaystyle \sigma^-\frac{d\tilde{q}^+}{dx_1} +\sigma^+\frac{d\tilde{q}^-}{dx_1}=\tilde{\mathfrak{g}}_1,\\[6pt]
\tilde{q}^+-\tilde{q}^-=\tilde{\mathfrak{g}}_2,\quad x_1=0,
\end{array}
\right. \label{10}
\end{align}
where $\tilde{\mathfrak{g}}_k$ is the Fourier transform of $\mathfrak{g}_k$ for $k=1,2$ and $\sigma^-$, $\sigma^+$ are the symbols of the operators $\Sigma^-_\gamma$, $\Sigma^+_\gamma$ respectively, that is
\[
\begin{split}
& \sigma^-=\sigma^-(\tau,\eta^\prime)=w_-^2-\widehat{\mathcal{E}}_1^2\eta^2+\varepsilon\big(\varepsilon\vert\widehat{\mathcal H}^\prime\vert^2\tau^2-2\widehat{\mathcal E}_1 i\tau w^{\perp}_-\big)=w_-^2-\widehat{\mathcal{E}}_1^2\eta^2+\mathcal{O}(\varepsilon)\,,\\
& \sigma^+=\sigma^+(\tau,\eta^\prime)=\ell^2+w_+^2,
\end{split}
\]
where
\[
\ell=\tau+i(\widehat{v}^\prime\cdot\eta^\prime)\,,\quad w_+=\widehat{H}'\cdot\eta',\quad w_-=\widehat{\mathcal{H}}'\cdot\eta'\,,\quad w^{\perp}_-=\widehat{\mathcal H}^\perp\cdot\eta^\prime\,.
\]
As usual, we define the hemi-sphere
\[
\Sigma:=\{(\tau,\eta^\prime)\in\mathbb C\times\mathbb R^2\,:\,\,\vert\tau\vert^2+\eta^2=1\,,\,\,\Re\tau\ge 0\}
\]
and denote by $\Xi$ the set of ``frequencies''
\[
\Xi:=\{(\tau, \eta^\prime)\in\mathbb C\times\mathbb R^2\,:\,\,\Re\tau\ge 0\,,\,\,(\tau,\eta^\prime)\neq (0,0)\}=]0,+\infty[\cdot\Sigma\,.
\]

Following \cite{MTT1}, we are going to construct a symbolic symmetrizer for the transformed problem \eqref{9}, \eqref{10}. We introduce the unknowns $\bold Y^\pm=\bold Y^\pm(\delta,x_1,\eta^\prime)$
\[
\bold Y^+=\begin{pmatrix} y_1^+ \\ y_2^+\end{pmatrix}=
\begin{pmatrix} \displaystyle \frac{d\tilde{q}^+}{dx_1} \\[3pt] \eta\tilde{q}^+\end{pmatrix}\quad \mbox{and}\quad
\bold Y^-=\begin{pmatrix} y_1^- \\ y_2^-\end{pmatrix}=
\begin{pmatrix} \displaystyle \frac{d\tilde{q}^-}{dx_1} \\[3pt] \sigma\tilde{q}^-\end{pmatrix},
\]
where $\sigma=\sigma(\tau,\eta^\prime)=\sqrt{\eta^2+\varepsilon^2\tau^2}$ denotes the {\em principal square root} of $\eta^2+\varepsilon^2\tau^2$, that is the square root of positive real part for $\Re\tau>0$, extended as a continuous function up to ``boundary frequencies'' $(\tau,\eta^\prime)\neq (0,0)$ with $\Re\tau=0$. Then, problem \eqref{9}, \eqref{10} is written as
\begin{align}
\frac{d}{dx_1}\bold Y=\mathcal{A}(\tau,\eta')\bold Y&\quad\mbox{for}\  x_1>0,\label{11}\\
\beta(\tau,\eta')\bold Y=\boldsymbol{\mathcal G} &\quad\mbox{at}\ x_1=0,\label{12}
\end{align}
where $\bold Y=(\bold Y^+,\bold Y^-)$, $\mathcal{A}={\rm diag}\,(\mathcal{A}^+,\mathcal{A}^-)$,
\begin{equation}\label{coeff_bound_data}
\mathcal{A}^+=\begin{pmatrix}0& \eta \\ \eta &0\end{pmatrix},\quad \mathcal{A}^-=\begin{pmatrix}0& \sigma \\ \sigma &0\end{pmatrix},\quad
\beta =\begin{pmatrix}\displaystyle\frac{\sigma^-}{\Lambda^2}& 0& \displaystyle\frac{\sigma^+}{\Lambda^2} &0 \\[9pt] 0&\displaystyle\frac{\sigma}{\Lambda} &0&-\displaystyle\frac{\eta}{\Lambda} \end{pmatrix},\quad \boldsymbol{\mathcal G}=\begin{pmatrix} \displaystyle\frac{\tilde{\mathfrak{g}}_1}{\Lambda^2} \\[9pt] \displaystyle\frac{\eta\sigma\tilde{\mathfrak{g}_2}}{\Lambda}\end{pmatrix}.
\end{equation}

\subsection{Lopatinski determinant}\label{LD}

The matrix $\mathcal{A}$ has the ``stable'' eigenvalues
\[
\lambda^+=-\eta\quad\mbox{and}\quad \lambda^-=-\sigma
\]
with the associated eigenvectors
\begin{equation}\label{stable_evc}
E^+ =\begin{pmatrix}1\\-1\\0\\ 0 \end{pmatrix} \quad\mbox{and}\quad E^-=\begin{pmatrix}0\\ 0 \\1\\-1\end{pmatrix}.
\end{equation}
Note that the matrix $\cal A(\tau,\eta^\prime)$ is diagonalizable for all $(\tau,\eta^\prime)\in\Xi$. More
precisely,
\[
 T {\cal A}(\tau, \eta^\prime) T^{-1} =
                                \begin{pmatrix} -\eta &0&0&0\\
                                                0&\eta &0&0\\
                                                0&0&-\sigma &0\\
                                                0&0&0&\sigma
                                                \end{pmatrix},\qquad
T=\frac{1}{2} \begin{pmatrix} 1&-1&0&0\\
                                1&1&0&0\\
                                0&0& 1&-1\\
0&0& 1& 1          \end{pmatrix} .
\]

Then, as in the hyperbolic theory, for problem \eqref{9}, \eqref{10} we define the Lopatinski determinant
\begin{equation}\label{delta}
\begin{split}
\Delta (\tau,\eta') &= \det \left[ \beta (E^+\ E^-)\right]=\det \begin{pmatrix}\displaystyle\frac{\sigma^-}{\Lambda^2}& \displaystyle\frac{\sigma^-}{\Lambda^2}\\[9pt] -\displaystyle\frac{\sigma}{\Lambda} &\displaystyle\frac{\eta}{\Lambda} \end{pmatrix}\\
&=\frac1{\Lambda^3}\left\{\eta \big(w_-^2-\widehat{\mathcal{E}}_1^2\eta^2+\varepsilon^2\vert\widehat{\mathcal H}^\prime\vert^2\tau^2-2\varepsilon\widehat{\mathcal E}_1i\tau w_-^{\perp}\big)+(\ell^2+w^2_+)\sigma \right\}\,.
\end{split}
\end{equation}
It is worthwhile noticing that $\beta(\tau,\eta^\prime)$ and $\Delta(\tau,\eta^\prime)$ defined above are homogeneous functions of degree zero with respect to $(\tau,\eta^\prime)$ in $\Xi$. Because of the homogeneity properties, one can reduce the study of the Lopatinski determinant to the hemi-sphere $\Sigma$, where it is a continuous function. If the Lopatinski determinant vanishes for $\Re\tau >0$, then the constant coefficients linearized  problem \eqref{27'}, \eqref{28'} is {\it ill-posed}, i.e. the piecewise constant basic state \eqref{constsol} is {\it unstable}. This never happens if the stability condition \eqref{stabcond} is satisfied, as it follows from \cite[Theorem 3.1]{T18}. Moreover, again from \cite{T18}, it can be proved the following proposition.

\begin{proposition}\label{prop}
Assume that \eqref{stabcond} holds. Then the equation $\Delta(\tau,\eta^\prime)=0$ has only simple roots $(\tau,\eta^\prime)\in\Sigma$ with $\Re\tau =0$.
\end{proposition}
\noindent
Arguing as in \cite{CS1} and \cite{MTT1}, we obtain the following result on the vanishing of the Lopatinski determinant.
\begin{lemma}\label{4.5}
Let $(\tau_0,\eta^\prime_0)\in\Sigma$ be a root of $\Delta(\tau,\eta^\prime)=0$. Then there exist a neighborhood $\mathcal{V}$
of $(\tau_0,\eta^\prime_0)$ in $\Sigma$ and a constant $k_0>0$ such that for all $(\tau,\eta^\prime)\in\mathcal{V}$ we have
$$
\left\vert\beta(\tau,\eta^\prime)(E^+,E^-)\bold Z\right\vert^2\geq
k_0\gamma^{2}\vert \bold Z\vert^2\qquad\forall\, \bold Z \in\mathbb{C}^2.
$$
\end{lemma}

Let us now state a technical result that will be used below in the construction of the symmetrizer.
\begin{proposition}\label{prop_sigma}
Let $\sigma=\sigma(\tau,\eta^\prime)=\sqrt{\eta^2+\varepsilon^2\tau^2}$ denote the principal square root of $\eta^2+\varepsilon^2\tau^2$ (that is the square root of positive real part for $\Re\tau>0$), extended as a continuous function up to boundary points $(\tau,\eta^\prime)\neq (0,0)$ with $\Re\tau=0$. Then
\begin{equation}\label{estimate_sigma}
\Re\sigma(\tau,\eta^\prime)\ge\frac{\varepsilon\gamma}{\sqrt 2}\,,\quad\forall\,(\tau,\eta^\prime)\in\Xi\,.
\end{equation}
\end{proposition}
\noindent
The proof will be given in \ref{proof_est_sigma}.

%
%\begin{proof}
%Let us rewrite the equation $\Delta(\tau,\boldsymbol{\omega})=0$ in terms of $s=i\tau \in\mathbb{C}$:
%\begin{equation}\label{eqs}
% (s-a^+)^2 + (s-a^-)^2 =( b^+)^2 + (b^-)^2
%\end{equation}
%($a^\pm, b^\pm$ are defined in \eqref{defabl}). Equation \eqref{eqs} is a quadratic equation for $s$ and has two roots
%\[
%s_{1,2}=\frac{1}{2}\left( (a^+ + a^-)^2 \pm \sqrt{ D(\boldsymbol{\omega})}\, \right),
%\]
%where $D(\boldsymbol{\omega})= 2((b^+)^2 + (b^-)^2) - (a^+ - a^-)^2$. Clearly, the equation $\Delta(\tau,\boldsymbol{\omega})=0$ has no unstable roots $\tau$ (of positive real part) if and only if
%both the roots $s_{1,2}$ are real, i.e., the quadratic form $D(\boldsymbol{\omega})$ (for $\omega_2$ and $\omega_3$) is
%nonnegative. This is so if and only if the stability condition \eqref{Syr1}, \eqref{Syr2} is satisfied.
%
%Under the sharpened stability condition \eqref{Syr3} the quadratic form $D(\boldsymbol{\omega})$ is positive definite
%(recall that \eqref{Syr3} implies \eqref{Syr1}) and the roots $s_{1,2}$ are distinct. These roots correspond to simple roots $\tau_{1,2}=-is_{1,2}$ of the equation $\Delta(\tau,\boldsymbol{\omega})=0$ with $\Re\tau=0$. That is, the uniform Lopatinski condition is violated.
%\end{proof}
%
\subsection{Construction of a degenerate symmetrizer}\label{symmetr}

This subsection will be entirely devoted to the construction of a symbolic symmetrizer of  \eqref{11}, \eqref{12}.
A general idea of symmetrizer for our elliptic-hyperbolic problem follows the same lines of the analogous construction made in \cite{MTT1}, which is inspired by the idea of Kreiss' symmetrizer \cite{Kreiss} for hyperbolic problems. We first reduce the ODE system in \eqref{11} to a diagonal form with the matrix $T\mathcal{A}T^{-1}$.
Then, multiplying the resulting system by a Herminian matrix $r(\tau,\eta^\prime)$ ({\it symmetrizer})
and using the boundary conditions and special properties of $r$, we derive the estimate
\begin{equation}\label{estY}
\vert \bold Y(\delta,0,\eta^\prime)\vert^2
\leq\frac{C}{\gamma^{2}} |\boldsymbol{\mathcal{G}}|^2 \Lambda^2
\end{equation}
by standard ``energy'' arguments. %Taking into account \eqref{sourcebound}, integrating estimate \eqref{estY} with respect to $(\delta,\eta^\prime)\in\mathbb{R}^3$, recalling the definition of $Y$, and using Plancherel's theorem, we obtain the desired estimate \eqref{estdotq}.

While constructing the symmetrizer we closely follow the plan and notation of Coulombel and Secchi
\cite{CS1}. The symbolic symmetrizer $r(\tau,\eta^\prime)$ of \eqref{11}, \eqref{12}
is sought to be a homogeneous function of degree zero with respect to $(\tau,\eta^\prime)\in\Xi$. Thus, it
is enough to construct $r(\tau,\eta^\prime)$ in the unit hemisphere $\Sigma$. Since the latter
is a compact set, by the use of a smooth partition of unity we still
reduce the construction of  $r(\tau,\eta^\prime)$ to that in a neighborhood of each point of $\Sigma$.
The analysis above (see subsection \ref{LD}) shows that we have to distinguish between three
different subclasses of frequencies $(\tau,\eta^\prime)\in\Sigma$ in the
construction of $r(\tau,\eta^\prime)$:
\begin{itemize}
\item [i.] The interior points $(\tau_0,\eta^\prime_0)$ of $\Sigma$ such
  that $\Re\tau_0>0$.
\item [ii.] The boundary points $(\tau_0,\eta^\prime_0)$ of $\Sigma$ where the Lopatinski
  condition is satisfied (i.e., $\Delta(\tau_0,\eta^\prime_0)\neq 0$).
\item[iii.] The boundary points $(\tau_0,\eta^\prime_0)$ where the Lopatinski
  condition breaks down (i.e., $\Delta(\tau_0,\eta^\prime_0)=0$).
\end{itemize}
The symmetrizer we are going to construct is degenerate in the sense that the uniform Lopatinski condition is violated and we have to treat case iii.

\subsubsection{Construction of the symmetrizer: the interior points
  (case i)}\label{kreiss1}

Let us consider a point $(\tau_0, \eta^\prime_0)\in \Sigma$ with $\Re \tau_0
>0$. Recall the matrix  ${\cal A}(\tau, \eta^\prime)$ is diagonalizable for
all $(\tau,\eta^\prime)\in\Xi$. In a neighborhood $\mathcal {V}$ of
$(\tau_0,\eta^\prime_0)$ the symmetrizer is defined by
\begin{equation}\label{symm_i}
r(\tau,\eta^\prime)=\begin{pmatrix}  -1&0&0&0 \\
                                    0& K&0&0\\
                                0&0&-1&0\\
                               0&0&0&K  \end{pmatrix} \qquad \forall\,
                                    (\tau,\eta^\prime)\in {\mathcal V},                                           \end{equation}
where $K\geq 1$ is a positive real number, to be fixed large enough. Let us set
$\Re M:=\frac{M+M^*}{2}$ for every complex square
matrix $M$. The matrix $r(\tau,\eta^\prime)$ is Hermitian and, in view of Proposition \ref{prop_sigma}, it satisfies
\begin{equation}\label{cond1symm}
\forall
(\tau,\eta^\prime)\in\mathcal{V},\quad\Re(r(\tau,\eta^\prime)T \mathcal{A}(\tau,\eta^\prime)T^{-1})\geq \kappa\varepsilon\eta I,
\end{equation}
 where $I$ denotes the identity matrix of order $4$ and where $\kappa$ is a suitable constant $0<\kappa\le 1$ for all $0<\varepsilon<\!<1$. As in \cite{MTT1}, the presence of $\eta$ in the right-hand side of inequality \eqref{cond1symm} must be understood as an ``elliptical degeneracy'' of the symmetrizer.
\newline
Furthermore, as in  \cite[Section 4.3.1]{MTT1}, for $K\geq 1$ sufficiently large the following inequality holds true
\begin{equation}\label{cond2symm}
\forall (\tau,\eta^\prime)\in\mathcal{V},\quad
r(\tau,\eta^\prime)+C\widetilde{\beta}^*(\tau,\eta^\prime)\widetilde{\beta}(\tau,\eta^\prime)\geq I,
\end{equation}
with a suitable positive constant $C$ and
$\widetilde{\beta}(\tau,\eta^\prime):=\beta(\tau,\eta^\prime)T^{-1}$
(we shrink the neighborhood $\mathcal{V}$ if necessary). We note that
the first and the third columns of the matrix $ T^{-1}$ are $E^+$ and $E^-$ in \eqref{stable_evc},
and the crucial point in obtaining inequality \eqref{cond2symm} is that the matrix $\beta(\tau,\eta^\prime)(E^+, E^-)$ is invertible because the Lopatinski determinant does not vanish
at $(\tau_0, \eta^\prime_0)$.
%%%%%%%%%%%%%%%%%%%%%%%%%%%%%%%%%%%%%%%%%%%%%%%%%%%%%%%%%%%%%%%%%%%%%%%%%%%%%%%%%%%%%%%%%%%
\subsubsection{The boundary points (case ii)}\label{kreiss2}
Let  $(\tau_0,\eta^\prime_0)$ belong to the subclass ii of $\Sigma$, namely,
$\Re\tau_0=0,$ and $\Delta(\tau_0,\eta^\prime_0)\neq 0$.
The symmetrizer
$r(\tau,\eta^\prime)$ is defined in a neighborhood of $(\tau_0,\eta^\prime_0)$ in a
completely similar manner as in case i, see \eqref{symm_i}. Similarly as in case i, one can prove that the symmetrizer satisfies the following inequalities:
\begin{equation}\label{cond1symm_a}
\forall
(\tau,\eta^\prime)\in\mathcal{V},\quad\Re(r(\tau,\eta^\prime)T \mathcal{A}(\tau,\eta^\prime)T^{-1})\geq \frac{\varepsilon}{\sqrt 2}\min\{\eta,\gamma\} I,
\end{equation}
\begin{equation}\label{cond2symm_b}
\forall (\tau,\eta^\prime)\in\mathcal{V},\quad
r(\tau,\eta^\prime)+C\widetilde{\beta}^*(\tau,\eta^\prime)\widetilde{\beta}(\tau,\eta^\prime)\geq I.
\end{equation}
with suitable constant $C>0$ and all $0<\varepsilon\ll 1$.
%%%%%%%%%%%%%%%%%%%%%%%%%%%%%%%%%%%%%%%%%%%%%%%%%%%%%%%%%%%%%%%%%%%%%%%%%%%%%%%%%%%%%%%%%%%%%%%
\subsubsection{The boundary points (case iii)}\label{kreiss3}
Let $(\tau_0,\eta^\prime_0)\in\Sigma$ be a point of type iii and denote
by $\mathcal {V}$ a neighborhood of $(\tau_0,\eta^\prime_0)$ in $\Sigma$. We define the symmetrizer in $\mathcal{V}$ by
$$
 r(\tau,\eta^\prime)=\begin{pmatrix}  -\gamma^2&0&0&0 \\
                                    0& K&0&0\\
                                0&0&-\gamma^2&0\\
                               0&0&0&K  \end{pmatrix} \qquad \forall\,
                                    (\tau,\eta^\prime)\in {\mathcal V},                                           $$
where $K\geq 1$ is a positive real number, to be fixed large enough.
 The matrix $r(\tau,\eta^\prime)$ above is Hermitian and  we have
\begin{equation}\label{Rr}
\begin{array}{l}
\Re(r(\tau,\eta^\prime)T\mathcal{A}(\tau,\eta^\prime)T^{-1})
\geq \displaystyle\frac{\varepsilon}{\sqrt 2}\min\{\eta,\gamma\}
\begin{pmatrix}
\gamma^{2} & 0 & 0 & 0\\
0 & 1 & 0 & 0\\
0 & 0 & \gamma^2 & 0\\
0 & 0 & 0 & 1
\end{pmatrix}.
\end{array}
\end{equation}

We also get that there exists a constant $C>0$ such that
\begin{equation}\label{eqC}
r(\tau,\eta^\prime)+C\widetilde{\beta}^{*}(\tau,\eta^\prime)\widetilde{\beta}(\tau,\eta^\prime)\geq\gamma^{2}I\qquad\forall\,(\tau,\eta^\prime)\in\mathcal{V}.
\end{equation}
The proof of \eqref{eqC} is based on Lemma \ref{4.5} concerning the vanishing of the Lopatinski determinant.
%%%%%%%%%%%%%%%%%%%%%%%%%%%%%%%%%%%%%%%%%%%%%%%%%%%%%%%%%%%%%%%%%%%%%%%%%%%%%%%%%%%%%%%%%%%%%%%%%
\subsection{Derivation of estimate \eqref{estY}}
We are now ready to derive estimate \eqref{estY}. Following \cite{CS1}, we introduce a smooth partition of
unity $\{\chi_j\}_{j=1}^J$ related to a given finite open covering
$\{\mathcal{V}_j\}_{j=1}^J$ of $\Sigma$. Namely, we have
\begin{equation*}
\begin{array}{l}
\chi_j\in C^{\infty},\quad
{\rm supp}(\chi_j)\subseteq\mathcal{V}_j,\ j=\overline{1,J},\quad \mbox{and}\quad
\sum\limits_{j=1}^J\chi^2_j\equiv 1.
\end{array}
\end{equation*}
Fix an arbitrary point $(\tau_0,\eta^\prime_0)\in\Sigma$ belonging to one of the classes
(i, ii or iii)  analyzed before and let $\mathcal{V}_j$ be an open
neighborhood of this point. We derive a local energy estimate in
$\mathcal{V}_j$ and then, by adding the resulting estimates over all
$j={\overline {1, J}}$, we obtain the desired global estimate.

${\bf 1^{\rm\bf st}}$ {\bf case.}
 Let $(\tau_0,\eta^\prime_0)$ belongs to class i or class ii.
We know from paragraphs \ref{kreiss1} and \ref{kreiss2} (see \eqref{cond1symm}, \eqref{cond2symm} and \eqref{cond1symm_a}) that there exist a $C^{\infty}$
mapping $r_j(\tau,\eta^\prime)$  defined on $\mathcal{V}_j$ such
that
\begin{itemize}
\item{} $r_j(\tau,\eta^\prime)$ is Hermitian;
\item{} the estimates
\begin{equation}\label{est1}
\begin{array}{l}
\Re\left(r_j(\tau,\eta^\prime)
  T \mathcal{A}(\tau, \eta^\prime)T^{-1}\right)\geq  K_j\varepsilon\min\{\eta,\gamma\} I,\\[3pt]
r_j(\tau,\eta^\prime)+C_j\widetilde{\beta}^*(\tau,\eta^\prime)\widetilde\beta(\tau,\eta^\prime)\geq I
\end{array}
\end{equation}
hold for all $(\tau,\eta^\prime)\in\mathcal{V}_j$, where $K_j$, $C_j$ are positive constants and we recall that $\widetilde{\beta}(\tau,\eta^\prime):=\beta(\tau,\eta^\prime)T^{-1}$ (the trivial estimate $\eta\ge\min\{\eta,\gamma\}$ is used in the right-hand side of \eqref{cond1symm}).
\end{itemize}
We set $\bold U_j(\tau,x_1,\eta^\prime):=\chi_j(\tau,\eta^\prime) T \bold Y(\delta,x_1,\eta^\prime)$.  Since $\chi_j$ is  supported on $\mathcal{V}_j$, we may think about $r_j$ extended by
zero to the whole of $\Sigma$. Then we extend $\chi_j$ and $r_j$ to the whole set of frequencies $\Xi$ as homogeneous
mappings of degree zero with respect to $(\tau,\eta^\prime)$.  Thus, from
equations \eqref{11}, \eqref{12} we obtain that $\bold U_j$ satisfies
\begin{equation}\label{systemUj}
\begin{cases}
{\displaystyle \frac{d\bold U_j}{dx_1}=T \mathcal{A}(\tau,\eta^\prime)T^{-1}\bold U_j,\quad x_1>0,}\\[3pt]
\widetilde\beta(\tau,\eta^\prime)\bold U_j(0)=\chi_j\,  \boldsymbol{\mathcal{G}}.
\end{cases}
\end{equation}
Taking the scalar product of the ODE system in \eqref{systemUj} with $r_j\bold U_j$,
integrating over $\mathbb{R}^+$ with respect to $x_1$,  and considering the
real part of the resulting equality, we are led to
\[
-\frac12(r_j(\tau,\eta^\prime)\bold U_j(\tau,0,\eta^\prime),\bold U_j(\tau,0,\eta^\prime))
=\int\limits_0^{+\infty}\Re\left(r_j(\tau,\eta^\prime)T
  \mathcal{A}(\tau,\eta^\prime)T^{-1} \bold U_j(\tau,x_1,\eta^\prime),\bold U_j(\tau,x_1,\eta^\prime)\right)dx_1.
\]
Then, by using estimates \eqref{est1} and the boundary condition in \eqref{systemUj}, one gets
\begin{equation*}
K_j\min\{\eta,\gamma\} \int\limits_0^{+\infty}\vert \bold U_j(\tau,x_1,\eta^\prime)\vert^2
dx_1+\frac12\vert
\bold U_j(\tau,0,\eta^\prime)\vert^2\leq\frac{C_j}{2}\chi_j^2(\tau,\eta^\prime)
|\boldsymbol{\mathcal {G}}|^2 .
\end{equation*}
Recalling the definition of $\bold U_j$, we obtain
\begin{equation}\label{4.32}
K_j\chi_j^2(\tau,\eta^\prime)\min\{\eta,\gamma\} \int\limits_0^{+\infty}\vert
\bold Y(\delta,x_1,\eta^\prime)\vert^2dx_1+\chi_j^2(\tau,\eta^\prime)\vert \bold Y(\delta,0,\eta^\prime)\vert^2
\leq C_j\chi_j^2(\tau,\eta^\prime)|\boldsymbol{\mathcal {G}}|^2.
\end{equation}

${\bf 2}^{\rm\bf  nd}$ {\bf case.}
It remains to prove a counterpart of estimate \eqref{4.32} for a neighborhood of a point $(\tau_0,\eta^\prime_0)\in \Sigma$ belonging to class ii, that is such that $\Re\,\tau_0=0$ and $\Delta(\tau_0,\eta^\prime_0)=0$. Let $\mathcal{V}_j$ be an open neighborhood of this
$(\tau_0,\eta^\prime_0)$ and $\chi_j$ the associated cut-off function. As was shown in paragraph \ref{kreiss3}
(see \eqref{Rr} and \eqref{eqC} and recall that $\gamma=\Re\tau\le 1$ for every $(\tau,\eta^\prime)\in\Sigma$), there exists a $C^{\infty}$
mapping $r_j(\tau,\eta^\prime)$ defined in $\mathcal{V}_j$, such that the following holds true
\begin{itemize}
\item{} $r_j(\tau,\eta^\prime)$ is Hermitian,
\item{} the estimates
\begin{equation}\label{est2}
\begin{array}{l}
\Re\left(r_j(\tau,\eta^\prime)T \mathcal{A}(\tau,\eta^\prime)T ^{-1}\right)\geq
\displaystyle\frac{\varepsilon}{\sqrt 2}\min\{\eta,\gamma\}\gamma^2 I,\\[6pt]
r_j(\tau,\eta^\prime)+C_j\widetilde{\beta}^*(\tau,\eta^\prime)\widetilde{\beta}(\tau,\eta^\prime)\geq \hat C_j\gamma^{2}I
\end{array}
\end{equation}
hold for all $(\tau,\eta^\prime)\in\mathcal{V}_j$, with some positive constants $C_j$, $\hat C_j$.
\end{itemize}
Recall that $r_j(\tau,\eta^\prime)$, $\mathcal{A}(\tau,\eta^\prime)$, and
$\beta(\tau,\eta^\prime)$ are assumed to be zero outside
$\mathcal{V}_j$. Then, we extend $r_j(\tau,\eta^\prime)$ and $\chi_j(\tau,\eta^\prime)$ to the whole of $\Xi$ as homogeneous mappings of degree $2$ and $0$ respectively. Hence, inequalities
\eqref{est2} become
\begin{equation}\label{fullest2}
\begin{array}{l}
\Re\left(r_j(\tau,\eta^\prime)T \mathcal{A}(\tau,\eta^\prime)T ^{-1}\right)\geq \displaystyle\frac{\varepsilon}{\sqrt 2}\min\{\eta,\gamma\}\gamma^2 I,\\[6pt]
r_j(\tau,\eta^\prime)+C_j\Lambda^2 \widetilde{\beta}^*(\tau,\eta^\prime)\widetilde{\beta}(\tau,\eta^\prime)\geq\hat C_j\gamma^{2}I
\end{array}
\end{equation}
for all $(\tau,\eta^\prime)\in\Xi$.

We again define $\bold U_j(\tau,x_1, \eta^\prime):=\chi_j(\tau,\eta^\prime)T \bold Y(\delta,x_1,\eta^\prime)$. Reasoning as above, we derive the estimate
\begin{equation}\label{4.35}
\displaystyle\frac{\varepsilon}{\sqrt 2}\min\{\eta,\gamma\}\chi^2_j(\tau,\eta^\prime) \int\limits_0^{+\infty}\vert
\bold Y(\delta,x_1,\eta^\prime)\vert^2dx_1+\hat{C}_j\chi^2_j(\tau,\eta^\prime)\vert \bold Y(\delta,0,\eta^\prime)\vert^2
\leq\frac{C_j}{\gamma^{2}}\,\chi^2_j(\tau,\eta^\prime)\Lambda^2(\tau,\eta^\prime)|\boldsymbol{\mathcal {G}}|^2,
\end{equation}
with a suitable positive constants $C_j$, $\hat C_j$.

\medskip
We now add up estimates \eqref{4.32} and \eqref{4.35} and use the fact that  $\{\chi_j\}$  is  a partition of unity. This leads us to the global estimate
\[
K\varepsilon \min\{\eta,\gamma\}\int\limits_0^{+\infty}\vert \bold Y (\delta,x_1,\eta^\prime)\vert^2dx_1+\hat C\vert \bold Y(\delta,0,\eta^\prime)\vert^2
\leq  C |\boldsymbol{\mathcal {G}}|^2 +\frac{C}{\gamma^{2}} \, \Lambda^2(\tau,\eta^\prime)|\boldsymbol{\mathcal {G}}|^2.
\]
Because of the inequality $ \Lambda(\tau,\eta^\prime)\geq\gamma$ we finally get
\[
K\varepsilon \min\{\eta,\gamma\}\int\limits_0^{+\infty}\vert \bold Y (\delta,x_1,\eta^\prime)\vert^2dx_1+\hat C\vert \bold Y(\delta,0,\eta^\prime)\vert^2 \leq\frac{C}{\gamma^{2}}\, \Lambda^2(\tau,\eta^\prime)|\boldsymbol{\mathcal {G}}|^2.
\]
The last estimate yields the desired estimate \eqref{estY}. To end up, we integrate \eqref{estY} with respect to $(\delta,\eta^\prime)$ on $\mathbb R^3$ to get
\[
\begin{split}
\int\limits_{\mathbb R^3}\vert \bold Y(\delta,0,\eta^\prime)\vert^2 d\delta d\eta^\prime & = \int\limits_{\mathbb R^3}\left\{\left\vert \frac{d\tilde{q}^+}{dx_1}\right\vert^2  + \left\vert \frac{d\tilde{q}^-}{dx_1}\right\vert^2 +\vert\eta\tilde{q}^+\vert^2+ \vert\sigma\tilde{q}^+\vert^2\right\} d\delta d\eta^\prime\\
&\leq\frac{C}{\gamma^{2}}\int\limits_{\mathbb R^3}\Lambda^2(\tau,\eta^\prime)|\boldsymbol{\mathcal {G}}(\delta,\eta^\prime)|^2 d\delta d\eta^\prime\,.
\end{split}
\]
On the other hand, in view of \eqref{coeff_bound_data} and Parseval's identity,
\[
\begin{split}
\int\limits_{\mathbb R^3}&\Lambda^2(\tau,\eta^\prime)|\boldsymbol{\mathcal {G}}(\delta,\eta^\prime)|^2 d\delta d\eta^\prime=\int\limits_{\mathbb R^3} \Lambda^2(\tau,\eta^\prime)\left\{\left\vert\frac{\tilde{\mathfrak{g}}_1}{\Lambda^2}\right\vert^2+\left\vert\displaystyle\frac{\eta\sigma\tilde{\mathfrak{g}}_2}{\Lambda}\right\vert^2\right\}  d\delta d\eta^\prime\\
&=\int\limits_{\mathbb R^3} \left\{\frac{\tilde{\vert\mathfrak{g}}_1\vert^2}{\Lambda^2}+\displaystyle\eta^2\vert\sigma\vert^2\vert\tilde{\mathfrak{g}}_2\vert^2\right\}d\delta d\eta^\prime
\le C\int\limits_{\mathbb R^3} \left\{\frac{\tilde{\vert\mathfrak{g}}_1\vert^2}{\Lambda^2}+\displaystyle\Lambda^4\vert\tilde{\mathfrak{g}}_2\vert^2\right\}d\delta d\eta^\prime\\
&=\Vert\mathfrak{g}_1\Vert^2_{-1,\gamma}+\Vert \mathfrak{g}_2\Vert^2_{2,\gamma}\,.
\end{split}
\]
From the above inequalities and again by Parseval's indentity, we deduce:
\begin{equation}\label{norm_der_estimate_0}
\Vert \nabla q^+\vert_{\omega}\Vert^2_{L^2(\omega)}+\Vert \partial_1 q^-\vert_{\omega}\Vert^2_{L^2(\omega)}\leq\frac{C}{\gamma^{2}}\{\Vert\mathfrak{g}_1\Vert^2_{-1,\gamma}+\Vert \mathfrak{g}_2\Vert^2_{2,\gamma}\}\,.
\end{equation}
Finally, using \eqref{trace}, \eqref{06} and the definition of $\mathfrak{g}_1$, $\mathfrak{g}_2$ see \eqref{gg}, we get
\[
\begin{split}
&\Vert\mathfrak{g}_1\Vert_{-1,\gamma}\le C\{\Vert f_{1,1}\vert_{\omega}\Vert_{1,\gamma}+\Vert\partial_1\breve{q}\vert_{\omega}\Vert_{1,\gamma}\}\le\frac{C}{\gamma}\nt f_{1,1}\nt_{2,\gamma}+C\nt f_1\nt_{2,\gamma}\le C\nt f_1\nt_{2,\gamma}\,,\\
&\Vert\mathfrak{g}_2\Vert_{2,\gamma}=\Vert\breve q\vert_{\omega}\Vert_{2,\gamma}\le C\nt f_1\nt_{3,\gamma}\,.
\end{split}
\]
Using the last inequalities to estimate the right-hand side of \eqref{norm_der_estimate_0}, we obtain
\begin{equation*}
\Vert \nabla q^+\vert_{\omega}\Vert^2_{L^2(\omega)}+\Vert \partial_1 q^-\vert_{\omega}\Vert^2_{L^2(\omega)}\leq\frac{C}{\gamma^{2}}\nt f_1\nt_{3,\gamma}^2
\end{equation*}
and adding \eqref{07}
\begin{equation}\label{gradqq-}
\Vert \nabla q\vert_{\omega}\Vert^2_{L^2(\omega)}+\Vert \partial_1 q^-\vert_{\omega}\Vert^2_{L^2(\omega)}\leq\frac{C}{\gamma^{2}}\nt f_1\nt_{3,\gamma}^2
\end{equation}

%{\bf DA QUI RIPRENDE LA VERSIONE DI YURI !!!}
%
%......................\\ ......................\\ ......................\\ ......................\\  ......................\\  ......................\\
%We finally deduce the estimate
%\begin{equation}\label{q+q-}
%\|\nabla q^+\|^2_{L^2(\omega )} + \|\nabla q^-\|^2_{L^2(\omega )} \leq \frac{C}{\gamma^2}\left(
%\|\mathfrak{g}_1\|^2_{2,\gamma} +\|\mathfrak{g}_2\|^2_{1,\gamma}\right).
%\end{equation}
%It follows from \eqref{06}, \eqref{07}, \eqref{gg} and the trace theorem that
%\begin{equation}\label{gradqq-}
%\|\nabla q\|^2_{L^2(\omega )} + \|\nabla q^-\|^2_{L^2(\omega )} \leq \frac{C}{\gamma^2}\nt f_1\nt^2_{4,\gamma}.
%\end{equation}

Restricting \eqref{U} to the boundary, by standard arguments we get the following estimate for the trace of U:
\begin{equation}\label{trU}
\gamma \|U|_{\omega}\|^2_{L^2(\omega )}\leq \frac{C}{\gamma} \left(\|\nabla q|_{\omega}\|^2_{L^2(\omega )}+\|f|_{\omega}\|^2_{L^2(\omega )}\right) .
\end{equation}
From \eqref{UV}, \eqref{phi6}, \eqref{gradqq-}, \eqref{trU} and the last two boundary conditions in \eqref{28''} we derive the estimate \eqref{est1.1} which implies \eqref{est0}. This completes the proof of Theorem \ref{mainth1}.

\section*{Acknowledgements}
A part of this work was done during a short stay  of Yuri Trakhinin  at the University of Brescia in March 2019. This author gratefully thanks  the Mathematical Division of the Department of Civil, Environmental, Architectural Engineering and Mathematics  of the University of Brescia for its kind hospitality. This work was also partially supported by RFBR (Russian Foundation for Basic Research) grant No. 19-01-00261-a and FFABR (Fondo di Finanziamento per le Attivit\`a Base di Ricerca) grant No. D83C18000060001.

\appendix
\section{Proof of Proposition \ref{prop_sigma}}\label{proof_est_sigma}
For fixed positive $\varepsilon$, direct calculations yield that the unique square root of $\varepsilon^2\tau^2+\eta^2$ with positive real part for $(\tau=\gamma+i\delta,\eta^\prime)\in\mathbb C\times\mathbb R^2$ such that $\gamma>0$ is
\begin{equation}\label{pos_real_part}
\sigma(\tau,\eta^\prime)=\sqrt{\frac{\alpha+r}{2}}+i\,\text{sgn}\,\beta\sqrt{\frac{-\alpha+r}{2}}\,,
\end{equation}
where
\begin{equation}\label{abr}
\begin{split}
& \alpha:=\Re(\varepsilon^2\tau^2+\eta^2)=\varepsilon^2(\gamma^2-\delta^2)+\eta^2\,,\quad\beta:=\Im(\varepsilon^2\tau^2+\eta^2)=2\varepsilon^2\gamma\delta\,,\\
& r:=\vert \varepsilon^2\tau^2+\eta^2\vert=\sqrt{\left(\varepsilon^2(\gamma^2-\delta^2)+\eta^2\right)^2+4\varepsilon^4\gamma^2\delta^2}
\end{split}
\end{equation}
and it is set
\[
\text{sgn}\,\theta=1\,,\quad\mbox{if}\,\,\theta\ge 0\qquad\mbox{and}\qquad\text{sgn}\,\theta=-1\,,\quad\mbox{if}\,\,\theta<0\,.
\]
Furthermore, the extension of $\sigma(\tau,\eta^\prime)$ to all boundary points $(i\delta,\eta^\prime)\in i\mathbb R\times\mathbb R^2$ such that $(\delta,\eta^\prime)\neq (0,0)$ is provided by
\begin{equation}\label{zero_real_part}
\sigma(i\delta,\eta^\prime)=\begin{cases}
\sqrt{-\varepsilon^2\delta^2+\eta^2}\,,\quad\mbox{if}\,\,-\varepsilon^2\delta^2+\eta^2\ge 0\,,\\
i\,\text{sgn}\,\delta\sqrt{\varepsilon^2\delta^2-\eta^2}\,,\quad\mbox{otherwise.}
\end{cases}
\end{equation}
From \eqref{zero_real_part}
\[
\Re\sigma(i\delta,\eta^\prime)\ge 0\,,\quad\forall\,(i\delta,\eta^\prime)\in i\mathbb R\times\mathbb R^2\,,\,\,(\delta,\eta^\prime)\neq (0,0)\,,
\]
that is \eqref{estimate_sigma} with $\gamma=0$. On the other hand, for $\tau=\gamma+i\delta$ with $\gamma>0$ and any $\eta^\prime\in\mathbb R^2$ we directly compute:
\[
\begin{split}
r^2 &=\left(\varepsilon^2(\gamma^2-\delta^2)+\eta^2\right)^2+4\varepsilon^4\gamma^2\delta^2=\varepsilon^4(\gamma^2-\delta^2)^2+2\varepsilon^2(\gamma^2-\delta^2)\eta^2+\eta^4+4\varepsilon^4\gamma^2\delta^2\\
&=\varepsilon^4\gamma^4+2\varepsilon^4\gamma^2\delta^2+\varepsilon^4\delta^4+2\varepsilon^2\gamma^2\eta^2-2\varepsilon^2\delta^2\eta^2+\eta^4=\varepsilon^4\gamma^4+2\varepsilon^4\gamma^2\delta^2+2\varepsilon^2\gamma^2\eta^2+(\varepsilon^2\delta^2-\eta^2)^2\\
&\ge (\varepsilon^2\delta^2-\eta^2)^2\,,
\end{split}
\]
hence, from \eqref{abr},
\[
\alpha+r\ge \varepsilon^2\gamma^2+\vert \varepsilon^2\delta^2-\eta^2\vert-(\varepsilon^2\delta^2-\eta^2)\ge\varepsilon^2\gamma^2
\]
and, from \eqref{pos_real_part},
\[
\Re\sigma(\tau,\eta^\prime)=\sqrt{\frac{\alpha+r}{2}}\ge\frac{\varepsilon}{\sqrt 2}\gamma\,,
\]
that is \eqref{estimate_sigma}.

\end{document}